\documentstyle[12pt]{article}
\def\Empty{}
\catcode`\@=11
\def\section{\@startsection {section}{1}{\z@}{-3.5ex plus-1ex minus
    -.2ex}{2.3ex plus.2ex}{\large\bf}}
\def\subsection{\@startsection{subsection}{2}{\z@}{-3.25ex plus-1ex
    minus-.2ex}{1.5ex plus.2ex}{\large\bf}}
\def\subsubsection{\@startsection{subsubsection}{3}{\z@}{-3.25ex plus
 -1ex minus-.2ex}{1.5ex plus.2ex}{\normalsize\bf}}

\def\eqalign#1{\,\vcenter{\openup\jot\m@th
  \ialign{\strut\hfil$\displaystyle{##}$&$\displaystyle{{}##}$\hfil
        \crcr#1\crcr}}\,}

\def\mydesc{\list{}{\labelwidth\z@ \itemindent-\leftmargin
\listparindent 1.5em
\let\makelabel\descriptionlabel}}

\def\fnum@figure{{\small Figure \thefigure}}
\def\fakefigure{\def\@captype{figure}}

\long\def\@makecaption#1#2{
    \vskip 10pt
    \def\FCap{#2} \def\NoCap{\ignorespaces}
    \ifx \FCap\NoCap
       \setbox\@tempboxa\hbox{#1}  % This is to avoid the damn colon.
      \else
       \setbox\@tempboxa\hbox{#1: \small \it #2}
    \fi
    \ifdim \wd\@tempboxa >\hsize   % IF longer than one line:
        \unhbox\@tempboxa\par      %   THEN set as ordinary paragraph.
      \else                        %   ELSE  center.
        \hbox to\hsize{\hfil\box\@tempboxa\hfil}
    \fi}

\def\@oddhead{\hbox{}\rightmark \hfil \rm\thepage}% Right heading.
\def\sectionmark#1{\markright {\sc{\ifnum \c@secnumdepth >\z@
      \S\thesection.\hskip 1em\relax \fi #1}}}

\catcode`\@=12

\def\oplabel#1{
  \def\OpArg{#1} \ifx \OpArg\Empty {} \else
        \label{#1}
  \fi}

\def\MakeStEnv#1{
  \newenvironment{#1}[2]{
  \begin{#1St} \oplabel{##1}%
  \global\def\CrntSt{\thetheoremSt}%
  {\rm ##2}%
}{
  \end{#1St} }
}
\MakeStEnv{theorem}
\MakeStEnv{corollary}
\MakeStEnv{proposition}
\MakeStEnv{lemma}
\MakeStEnv{conjecture}

\newenvironment{proof}[1]{
  \def\PfArg{#1}
  \ifx\PfArg\Empty
        \edef\PfArg{\CrntSt}  \fi
 \startproof{\PfArg}%
}{
  \finishproof{\PfArg}
}
\newcommand{\startproof}[1]{
  \medbreak\mbox{}
  {\it Proof of #1:}%
}
\newcommand{\finishproof}[1]{
  \def\FPArg{#1}
  \ifx\FPArg\Empty
        \def\FPArg{\CrntSt}  \fi
  \smallbreak\noindent\makebox[\textwidth]{\hfill\fbox{\FPArg}}
  \medbreak\noindent
}

\def\ol#1{\overline{#1}}     % overline

\def\eq{\!=\!}               % improved spacing for equals sign

\def\norm#1{\|#1\|}    % norm  || #1 ||
\def\rel{\hbox{\ rel\ }}
\def\set#1{\{ #1 \}}         % set braces  {  }
\def\vbar{\hbox{\hbox{$\;\vert\;$}}}% vertical bar with wider
\def\Isom{\hbox{\it Isom\/}}

\def\SO{\hbox{SO}}

\def\today{\ifcase\month\or
   January\or February\or March\or April\or May\or June\or
   July\or August\or September\or October\or November\or December\fi
   \space\number\day, \number\year}

\newcommand{\hyperbolic}{{\bf H}}

\newcommand{\integers}{{\bf Z}}

\newcommand{\reals}{{\bf R}}

\def\Z{\integers}
\def\R{\reals}

\def\H{\hyperbolic}

\def\Diff{\hbox{\it Diff\/}}
\def\diff{\hbox{\it diff\/}}
\def\Fr{\hbox{\it Fr\/}}
\def\Imb{\hbox{\it Imb\/}}
\def\imb{\hbox{\it imb\/}}

\def\Exp{\hbox{\rm Exp}}
\def\Maps{\hbox{\rm Maps}}

\def\mapdown#1{\big\downarrow % makes a downarrow with a mathematics
            \rlap            % mode symbol to the right
            {\smash{$\vcenter       % (for use in the "diagram" template)
            {\hbox{$         %
            \scriptstyle#1   %
            $}}$}}}           %

\def\mapright#1{\smash{      % makes a long rightarrow with a mathematics
            \mathop          % mode symbol above it
            {\longrightarrow % (for use in the "diagram" template)
            }\limits^{#1}}}  %

\newcommand{\marginwrite}[1]{}

\begin{document}

\title{Fiber-preserving imbeddings and diffeomorphisms}
\author{John Kalliongis and Darryl McCullough}
\date{{\footnotesize Department of Mathematics, Saint Louis
University, St.~Louis, MO 63103}
\\
{\footnotesize kalliongisje@sluvca.slu.edu}
\ \\
\ \\
{\footnotesize Department of Mathematics, University of Oklahoma,
Norman, OK 73019}
{\footnotesize dmccullough@math.ou.edu}
\ \\
\ \\
{\footnotesize \today
}}
\maketitle
\vfill
\noindent {\footnotesize 1991 Mathematics Subject Classification\ \
Primary: 57R35\ \ Secondary: 57M99}

\section{Introduction}
\label{intro}
\marginwrite{intro}

Let $\Diff(M)$ be the group of diffeomorphisms of a smooth manifold
$M$, with the $C^\infty$-topology. For a smooth submanifold $N$ of
$M$, denote by $\Imb(N,M)$ the space of all smooth imbeddings $j$ of
$N$ into $M$ such that $j^{-1}(\partial M)\eq N\cap \partial M$.  In
\cite{P}, R.~Palais proved a useful result relating diffeomorphisms
and imbeddings. In the case when $M$ is closed, it says that if
$W\subset V$ are submanifolds of $M$, then the mappings $\Diff(M)\to
\Imb(V,M)$ and $\Imb(V,M)\to \Imb(W,M)$ obtained by restricting
diffeomorphisms and imbeddings are locally trivial, and hence are
Serre fibrations. The same results, with variants for manifolds with
boundary and more complicated additional boundary structure, were
proven by J.~Cerf in~\cite{Cerf}. Among various applications of these
results, the Isotopy Extension Theorem follows by lifting a path in
$\Imb(V,M)$ starting at the inclusion map of $V$ to a path in
$\Diff(M)$ starting at $1_M$. Moreover, parameterized versions of
isotopy extension follow just as easily from the homotopy lifting
property for $\Diff(M)\to \Imb(V,M)$ (see corollary~\ref{isotopy
lifting}).

In the common situation of a fibering of manifolds, it is natural to
consider the spaces of imbeddings and diffeomorphisms that respect the
fibered structure. Consider a (smooth) fibering $p\colon E\to B$ of
compact manifolds, possibly with boundary. (Actually, most of our
results allow $E$ and $B$ to be noncompact, although the fiber and the
relevant submanifolds must be assumed to be compact. Also, we prove
versions with control relative to subsets of the boundary of $B$ and
their preimages in $E$. For clarity we omit such complications in
this introductory discussion.) A diffeomorphism of $E$ is called {\it
fiber-preserving} when it takes each fiber of $E$ to a fiber of $E$,
and {\it vertical} when it takes each fiber to itself. The space
$\Diff_f(E)$ of fiber-preserving diffeomorphisms of $E$ contains the
subspace $\Diff_v(E)$ of vertical diffeomorphisms. Any
fiber-preserving diffeomorphism $g$ of $E$ induces a diffeomorphism
$\ol{g}$ of $B$, and this defines a map from $\Diff_f(E)$ to
$\Diff(B)$ for which the preimage of the identity map is $\Diff_v(E)$.
In section~\ref{project} we prove

\medskip
\noindent{\bf Projection Theorem} (Theorem \ref{project diffs}) {\em
$\Diff_f(E)\to\Diff(B)$ is locally trivial. }
\medskip

\noindent This theorem is essentially due to W. Neumann and F. Raymond
(see the comments below). The homotopy extension property for the
projection fibration translates directly into the following.

\medskip
\noindent{\bf Parameterized Isotopy Extension Theorem} (Corollary
\ref{isotopy lifting}) {\em Suppose that $p\colon E\to B$ is a
fibering of compact manifolds, and suppose that for each $t$ in a
path-connected parameter space $P$, there is an isotopy $g_{t,s}$
such that $g_{t,0}$ lifts to a diffeomorphism $G_{t,0}$ of $E$. Assume
that sending $(t,s)\to g_{t,s}$ defines a continous function from
$P\times [0,1]$ to $\Diff(B)$ and sending $t$ to $G_{t,0}$ defines a
continuous function from $P$ to $\Diff(E)$. Then the family $G_{t,0}$
extends to a continuous family on $P\times I$ such that for each
$(t,s)$, $G_{t,s}$ is a fiber-preserving diffeomorphism inducing
$g_{t,s}$ on~$B$.}
\medskip

\noindent
A submanifold of $E$ is called {\it vertical} if it is a union of
fibers, and in this case it will be assumed to have the fibered
structure so that the inclusion map is fiber-preserving. An imbedding
of a fibered manifold $W$ into $E$ is called {\it fiber-preserving}
when the image of each fiber of $W$ is a fiber of $E$. The space of
all fiber-preserving imbeddings from $W$ to $E$ is denoted by
$\Imb_f(W,E)$. When $W\subseteq E$, $\Imb_f(W,E)$ contains the
subspace of {\it vertical} imbeddings $\Imb_v(W,E)$ which take each
fiber to itself. For fiber-preserving and vertical imbeddings of
vertical submanifolds, we have a more direct analogue of Palais'
results.

\medskip
\noindent{\bf Restriction Theorem} (Corollary \ref{corollary3}) {\em
Let $V$ and $W$ be vertical submanifolds of $E$ with $W\subseteq V$,
each of which is either properly imbedded or codimension-zero. Then
the restrictions $\Imb_f(V,E)\to \Imb_f(W,E)$ and $\Imb_v(V,E)\to
\Imb_v(W,E)$ are locally trivial.}
\medskip

\noindent As shown in theorem \ref{square}, the Projection and
Restriction Theorems can be combined into a single commutative square
in which all four maps are locally trivial:
$$\vbox{\halign{\hfil#\hfil\quad&#&\quad\hfil#\hfil\cr
$\Diff_f(E)$&$\longrightarrow$&$\Imb_f(W,E)$\cr
\noalign{\smallskip}
$\mapdown{}$&&$\mapdown{}$\cr
\noalign{\smallskip}
$\Diff(B)$&$\longrightarrow$&$\Imb(p(W),B)\rlap{\ .}$\cr}}$$

In 3-dimensional topology, a key role is played by manifolds
admitting a more general kind of fibered structure, called a Seifert
fibering. Some general references for Seifert-fibered 3-manifolds are
\cite{Hempel,Jaco1,JS,Orlik,OVZ,Scott,Seifert,Waldhausen1,Waldhausen2}.
In section~\ref{sfiber}, we prove the analogues of the results
discussed above for many Seifert fiberings $p\colon\Sigma\to{\cal O}$,
not necessarily 3-dimensional. Actually, we work in a somewhat more
general context, called {\t singular fiberings}, which resemble
Seifert fiberings but for which none of the usual structure of the
fiber as a homogeneous space is required.

In the late 1970's fibration results akin to our Projection Theorem
for the singular fibered case were proven by W.~Neumann and
F.~Raymond~\cite{N-R}. They were interested in the case when $\Sigma$
admits an action of the $k$-torus~$T^k$ and $\Sigma\to{\cal O}$ is the
quotient map to the orbit space of the action. They proved that the
space of (weakly) $T^k$-equivariant homeomorphisms of $\Sigma$ fibers
over the space of homeomorphisms of ${\cal O}$ that respect the orbit
types associated to the points of ${\cal O}$. A detailed proof of this
result when the dimension of $\Sigma$ is $k+2$ appears in the
dissertation of C.~Park~\cite{Park}. Park also proves analogous
results for space of weakly $G$-equivariant maps for principal
$G$-bundles and for Seifert fiberings of arbitrary
dimension~\cite{Park,Park1}. These results do not directly overlap
ours since we always consider the full group of fiber-preserving
diffeomorphisms without any restriction to $G$-equivariant maps
(indeed, no assumption of a $G$-action is even present).

Some technical applications of our results appear in \cite{M-R}. In
the present paper we give one main application.  For a Seifert-fibered
manifold $\Sigma$, $\Diff(\Sigma)$ acts on the set of Seifert
fiberings, and the stabilizer of a given fibering is
$\Diff_f(\Sigma)$, thus the space of cosets
$\Diff(\Sigma)/\Diff_f(\Sigma)$ is the {\it space of Seifert
fiberings} of $\Sigma$. We prove in section~\ref{sfspace} that for a
Seifert-fibered Haken 3-manifold, each component of the space of
Seifert fiberings is weakly contractible (apart from a small list of
well-known exceptions, the space of Seifert fiberings is connected).
This result is originally due to Neumann and Raymond, since it is an
immediate consequence of the results in \cite{N-R} combined with
contemporaneous work of Hatcher~\cite{Hatcher}. We make the same use
of~\cite{Hatcher}.

Our results will be proven by adapting the method developed
in~\cite{P}. The main new idea needed for the fibered case is a
modification of the usual exponential map, called the {\it aligned}
exponential map $\Exp_a$. This is defined and discussed in
section~\ref{exponent}. Section~\ref{metrics} contains some
preliminaries needed for carrying out Palais' approach for manifolds
with boundary. In section~\ref{palais}, we reprove the main result of
\cite{P} for manifolds which may have boundary. This
duplicates~\cite{Cerf} (in fact, the boundary control there is more
refined than ours), but is included to furnish lemmas as well as to
exhibit a prototype for the approach we use to deal with the bounded
case in our later settings. In section~\ref{orbifold}, we give the
analogues of the results of Palais and Cerf for smooth orbifolds,
which for us are quotients $\widetilde{\cal O}/H$ where
$\widetilde{\cal O}$ is a manifold and $H$ is a group acting smoothly
and properly discontinuously on $\widetilde{\cal O}$. Besides being of
independent interest, these analogues are needed for the case of
singular fiberings.

By a submanifold $N$ of $M$, we mean a smooth submanifold. When $M$
has boundary and $\dim(N)<\dim(M)$, we always require that $N$ be
properly imbedded in the sense that $N\cap \partial M=\partial N$. If
$N$ has codimension~0, we require that the frontier of $N$ be a
codimension-1 submanifold of $M$. In particular, it is understood that
the elements of $\Imb(N,M)$ carry $N$ to a submanifold satisfying
these conditions. The notation $\Diff(M\rel\partial M)$ means the
space of diffeomorphisms which restrict to the identity map on each
point of $\partial M$, and for $X\subseteq M$, $\Imb(X,M\rel\partial
M)$ means the imbeddings that equal the inclusion on $X\cap\partial
M$. For $K\subseteq M$, $\Diff^K(M)$ means the diffeomorphisms that
agree with the identity on $M-K$. We say that $K$ is a neighborhood of
the subset $X$ when $X$ is contained in the topological interior of
$K$. If $K$ is a neighborhood of a submanifold $N$, then $\Imb^K(N,M)$
means the elements $j$ in $\Imb(N,M)$ such that $K$ is a neighborhood
of~$j(N)$.

The second author thanks the MSRI for its support while the present
manuscript was in preparation. Both authors appreciate the continued
support of St.~Louis University for their collaborative work.

\section{Metrics which are products near the boundary}
\label{metrics}
\marginwrite{metrics}

When $M$ has a Riemannian metric, we denote by $d$ the associated
topological metric defined by putting $d(x,y)$ equal to the infimum of
the lengths of all piecewise differentiable paths from $x$ to $y$ when
$x$ and $y$ lie in the same component of $M$, and equal to~1 if $x$
and $y$ lie in different components.

Let $V$ be a (possibly empty) compact submanifold of $M$. Recall that
we always assume that $V$ is properly imbedded, if it has positive
codimension, or that the frontier of $V$ is a properly imbedded
codimension-1 submanifold, if $V$ has codimension~0. Fix a smooth
collar $\partial M\times [0,2]$ of $\partial M$ such that $V\cap
\partial M\times [0,2]$ is a union of $[0,2]$-fibers. Such a collar
can be obtained by constructing an inward-pointing vector field on a
neighborhood of $\partial M$ which is tangent to $V$, and using the
integral curves associated to the vector field to produce the collar.
On $\partial M\times [0,2)$, fix a Riemannian metric that is the
product of a metric on $\partial M$ and the usual metric on $[0,2)$.
Form a metric on $M$ from this metric and any metric defined on all of
$M$ using a partition of unity subordinate to the open cover
$\set{\partial M\times[0,2), M-\partial M\times I}$, where
$I\eq[0,1]$. Such a metric is said to be a {\it product near $\partial
M$} such that $V$ {\it meets the collar $\partial M\times I$ in
$I$-fibers}. It has the following properties for $0\leq t\leq 1$:

\begin{enumerate}
\item[{\rm(i)}] If $x\in M-\partial M\times I$, then $d(x,\partial
M)>1$.
\item[{\rm(ii)}] If $x\eq (y,t)\in\partial M\times I$, then
$d(x,\partial M)\eq t$.
\item[{\rm(iii)}] If $V$ has positive codimension, then for any
tubular neighborhood of $V$ obtained by exponentiating a normal bundle
of $V$, the fiber at each point of $V\cap \partial M\times \set{t}$
lies in $\partial M\times\set{t}$. In particular, the fiber of each
point in $V\cap \partial M$ lies in~$\partial M$. If $V$ has
codimension $0$, then the corresponding statement holds for any
tubular neighborhood of the frontier of~$V$.
\end{enumerate}

A Riemannian metric is called {\it complete} if every Cauchy sequence
converges. For a complete Riemannian metric on $M$, a geodesic can be
extended indefinitely unless it reaches a point in the boundary of
$M$, where it may continue or it may fail to be extendible because it
``runs out of the manifold.''

One may obtain a complete metric on $M$ that is a product near
$\partial M$ such that $V$ meets the collar $\partial M\times I$ in
$I$-fibers as follows. Carry out the previous construction using a
metric on $\partial M\times[0,2)$ that is the product of a complete
metric on $\partial M$ and the standard metric on $[0,2)$. Define
$f\colon M-\partial M\to (0,\infty)$ by putting $f(x)$ equal to the
supremum of the values of $r$ such that $\Exp$ is defined on all
vectors in $T_x(M)$ of length less than~$r$. Let $g\colon M-\partial
M\to (0,\infty)$ be a smooth map that is an $\epsilon$-approximation
to $1/f$, and let $\phi\colon M\to[0,1]$ be a smooth map which is
equal to~$0$ on $\partial M\times I$ and is~$1$ on $M-\partial
M\times[0,2)$. Give $M\times[0,\infty)$ the product metric, and define
a smooth imbedding $i\colon M\to M\times[0,\infty)$ by $i(x)\eq
(x,\phi(x)g(x))$ if $x\notin\partial M$ and $i(x)\eq (x,0)$ if $x\in
\partial M$. The restricted metric on $i(M)$ agrees with the original
metric on $\partial M\times I$ and is complete.

From now on, all metrics will be assumed to be complete.

\section{The Palais-Cerf restriction theorem}
\label{palais}\marginwrite{palais}

In this section we modify some results from \cite{P} to apply to the
bounded case. This duplicates~\cite{Cerf}, in fact our results are not
as general since we do not work in the setting of manifolds with
corners. On the other hand, our argument will provide lemmas needed
for the fibered cases, and is the prototype for the approach we use to
deal with the bounded case in our later settings.

Let $X$ be a $G$-space and $x_0\in X$. A {\it local cross-section} for
$X$ at $x_0$ is a map $\chi$ from a neighborhood $U$ of $x_0$ into $G$
such that $\chi(u)x_0\eq u$ for all $u\in U$. By replacing $\chi(u)$
by $\chi(u)\chi(x_0)^{-1}$, one may always assume that $\chi(x_0)\eq
1_G$, If $X$ admits a local cross-section at each point, it is said to
admit local cross-sections. From \cite{P} we have

\begin{proposition}{theoremA}{} Let $G$ be a topological group and $X$ a
$G$-space admitting local cross-sections. Then any equivariant map of
a $G$-space into $X$ is locally trivial.
\marginwrite{theoremA}
\end{proposition}

\noindent In fact, when $\pi\colon Y\to X$ is $G$-equivariant, the
local coordinates on $\pi^{-1}(U)$ are just given by sending the point
$(u,z)\in U\times \pi^{-1}(y_0)$ to $\chi(u)\cdot z$. Some additional
properties of the bundles obtained in proposition~\ref{theoremA} are
given in~\cite{P}.

The following technical lemma will simplify some of our applications
of proposition~\ref{theoremA}.

\begin{proposition}{inclusion}{} Let $V$ be a submanifold of $M$, let
$I(V,M)$ be a space of imbeddings of $V$ into $M$, and let $D(M)$ a
group of diffeomorphisms of $M$. Suppose that for every $i\in I(V,M)$,
the space of imbeddings $I(i(V),M)$ has a local $D(M)$ cross-section
at the inclusion map of $i(V)$ into $M$. Then $I(V,M)$ has local
cross-sections.
\marginwrite{inclusion}
\end{proposition}

\begin{proof}{} Denote by $j_{i(V)}$ the inclusion map of $i(V)$ into
$M$. Let $i\in I(V,M)$ and define $Y\colon I(V,M)\to I(i(V),M)$ by
$Y(j)=ji^{-1}$. For a local cross-section $\chi\colon U\to D(M)$ at
$j_{i(V)}$, define $Y_1$ to be the restriction of $Y$ to $Y^{-1}(U)$,
a neighborhood of $i$ in $I(V,M)$. Then $\chi Y_1\colon Y^{-1}(U)\to
D(M)$ is a local cross-section for $I(V,M)$ at~$i$. For if $j\in
Y^{-1}(U)$ and $x\in V$, then $\chi (Y_1(j)) (i(x))= \chi
(Y_1(j))(j_{i(V)}(i(x))) = Y_1(j)(i(x))=j(x)$.
\end{proof}

The results in \cite{P} depend in large part on three lemmas, called
lemmas~b, c and~d there. Here, we adapt their statements and proofs to
the context of manifolds with boundary. First, for $L\subseteq M$ define
$\Maps^L(M,M)$ be the space of smooth maps $f\colon M\to M$ such that
$f(\partial M)\subseteq \partial M$ and $f(x)\eq x$ for all $x\in
M-L$.

\begin{lemma}{J}{(Palais' lemma b)} Let $K$ be a compact subset of a
manifold $M$. Then there exists a neighborhood $J$ of $1_M$ in
$\Maps^K(M,M)$ which consists of diffeomorphisms.
\marginwrite{J}
\end{lemma}

\begin{proof}{} There exists a ($C^\infty$-) neighborhood $N$ of the
identity consisting of maps $f$ for which the differential $T_x(f)$ is
an isomorphism for all $x$. Since $f(\partial M)\subseteq\partial M$,
this implies that $f$ is a local diffeomorphism. Since $K$ is compact,
the preimage of any compact subset of $M$ under $f$ is compact.
Therefore $f$ is a covering map. If $M\neq K$, then this covering must
be 1-fold off of $f(K)$, hence must be a diffeomorphism, so assume
that $M$ is compact. Fix $\epsilon>0$ such that no closed
noncontractible loop in $M$ has length less than $4\epsilon$, and let
$J$ consist of the elements in $N$ such that $d(f(x),x)<\epsilon$ for
all $x$. Suppose for contradiction that $f\in J$ but $f(p)\eq f(q)$
for $p\neq q$. Then $d(p,q)<2\epsilon$, and if $\alpha$ is a geodesic
from $p$ to $q$ of length less than $2\epsilon$, then the diameter of
$f(\alpha)$ is less than $4\epsilon$. Therefore $f(\alpha)$ is a
contractible loop, a contradiction.
\end{proof}

For the next lemmas, we set some notation. The projection from the
tangent bundle $T(M)$ to $M$ is denoted by~$\pi$. For a submanifold
$V$ of $M$, let ${\cal X}(V,T(M))$ denote the sections $X$ from $V$ to
$T(M)\vert_V$ such that
\begin{enumerate}
\item[(1)] if $x\in V\cap\partial M$, then $X(x)$ is tangent to
$\partial M$, and
\item[(2)] $\Exp(X(x))$ is defined for all $x\in V$.
\end{enumerate}

\noindent When the metric is a product near the boundary,
property~(1) implies that if $x\in V\cap \partial M$, then
$\Exp(X(x))\in\partial M$.

A zero section will usually be denoted by $Z$.  The vector fields
satisfying~(1) and~(2) (i.~e.~the case $V\eq M$) are denoted simply by
${\cal X}(T(M))$. When $L$ is a subset of $M$, denote by ${\cal
X}^L(T(M))$ the elements of ${\cal X}(T(M))$ which agree with $Z$
outside of~$L$. A subscript ``$<\delta$'' indicates the sections such
that each image vector has length less than the positive number
$\delta$, thus for example
$${\cal X}_{<1/2}(V,T(M))\eq
\set{X\in{\cal X}(V,T(M))\vbar \norm{X(x)}<1/2\hbox{\ for all
$x\in V$}}\ .$$

\begin{lemma}{extension}{(Palais' lemma c)} Let $V$ be a compact
submanifold of the smooth manifold $M$ and $L$ a neighborhood of $V$
in $M$. Assume that the metric on $M$ is a product near $\partial M$
such that $V$ is vertical in $\partial M\times I$. Then there exists a
continuous map $k\colon {\cal X}_{<1/2}(V,T(M))\to {\cal X}^L(T(M))$
such that $k(X)(x)\eq X(x)$ for all $x$ in $V$ and all $X\in {\cal
X}_{<1/2}(V,T(M))$. Moreover, $k(Z)\eq Z$, and if $S\subseteq\partial
M$ is a closed neighborhood in $\partial M$ of $S\cap \partial V$, and
$X(x)\eq Z(x)$ for all $x\in S\cap \partial V$, then $k(X)(x)\eq Z(x)$
for all $x\in S$.
\marginwrite{extension}
\end{lemma}

\begin{proof}{} Suppose first that $V$ has positive codimension. Let
$\nu_\epsilon(V)$ denote the vectors of length less than $\epsilon$ in
the normal bundle of $V$. Fix $\epsilon<1/2$ and sufficiently small so
that $j\colon\nu_\epsilon(V)\to M$ defined by exponentiation is a
tubular neighborhood of $V$ contained in $L$, and such that the union
of the fibers at points in $S\cap\partial V$ is contained in $S$, and
the union of the fibers at points in $(\partial M-S)\cap\partial V$ is
contained in $\partial M-S$.  By property (iii) of the metric, the
fiber of this neighborhood at each point of $V\cap \partial M\times
\set{t}$ lies in $\partial M\times\set{t}$. Since $V$ is compact, we
may choose $\epsilon$ sufficiently small so that $j(\omega)\in\partial
M\times I$ only when $\pi(\omega)\in\partial M\times I$.

Suppose $v\in T_x(M)$ and that $\Exp(v)$ is defined. For all $u\in
T_x(M)$ define $P(u,v)$ to be the vector that results from parallel
translation of $u$ along the path that sends $t$ to~$\Exp(tv)$, $0\leq
t\leq 1$. Note that $P(u,Z(x))\eq u$ for all $u$.  Let $\alpha\colon
M\to [0,1]$ be a smooth function which is identically~1 on $V$ and
identically~0 on $M-j(\nu_{\epsilon/2}(V))$. Define $k\colon
{\cal X}_{<1/2}(V,T(M))\to{\cal X}^L(T(M))$ by
$$k(X)(x)\eq
\cases{
\alpha(x)P(X(\pi(j^{-1}(x))),j^{-1}(x))&for $x\in j(\nu_{\epsilon}(V))$\cr
Z(x)&for $x\in M-j(\nu_{\epsilon/2}(V))$\ .\cr}$$

\noindent For $x\in V$, $j^{-1}(x)\eq Z(x)$ and $\alpha(x)\eq 1$,
so $k(X)(x)\eq X(x)$. We must check that $k(X)\in{\cal X}^L(T(M))$.
Since the metric on $M$ is a product near $\partial M$,
$k(X)$ satisfies condition~(1) to be in ${\cal X}(T(M))$. To verify
that it satisfies condition~(2), fix $x$ such that $k(X)(x)\neq
Z(x)$. Suppose first that $x\eq (y,t)\in\partial M\times I$. Then
$\pi(j^{-1}(x))$ has the form $(y',t)$. If $t\geq 1/2$ then
$\norm{P(X(\pi(j^{-1}(x))),j^{-1}(x))}\eq\norm{X(\pi(j^{-1}(x)))}<1/2$,
so since $\alpha(x)\leq 1$, $\norm{k(X)(x)} \leq t\eq d(x,\partial M)$
and $\Exp(k(X)(x))$ is defined.  Suppose $t\leq 1/2$. Since the metric
is a product near $\partial M$, the component of $k(X)(x)$ in the
$I$-direction can be identified with the component of
$X(\pi(j^{-1}(x)))$ in the $I$-direction, so $\Exp(k(X)(x))$ is
defined when $\Exp(X(\pi(j^{-1}(x))))$ is. Finally, if $x\notin
\partial M\times I$ then also $\pi(j^{-1}(x))\notin \partial M\times
I$ so $d(\pi(j^{-1}(x)),\partial M)>1$. Since
$\norm{X(\pi(j^{-1}(x)))}<1/2$ and $\norm{j^{-1}(x)}<\epsilon<1/2$,
$\Exp(k(X)(x))$ is a point that lies within distance~1
of~$\pi(j^{-1}(x))$. The fact that $k(Z)\eq Z$ is immediate from the
definition.  Finally, if $X(x)\eq Z(x)$ for all $x\in S\cap\partial
V$, then since the metric is a product near the boundary it follows
that $k(X)(x)\eq Z(x)$ for all $x\in S$.

Now suppose that $V$ has codimension zero, so that its frontier $W$ is
a properly imbedded submanifold. Let $\nu^+_\epsilon(W)$ denote the
vectors of lengths less than $\epsilon$ in the normal bundle of $W$
that exponentiate into $\overline{M-V}$. Proceed as before, but define
$k$ by
$$k(X)(x)\eq
\cases{X(x)&for $x\in V$\cr
\alpha(x)P(X(\pi(j^{-1}(x))),j^{-1}(x))&for $x\in j(\nu^+_{\epsilon}(V))$\cr
Z(x)&for $x\in M-j(\nu^+_{\epsilon/2}(V))$\ .\cr}$$
\end{proof}

For our proof of lemma~d, we introduce some additional notation.
Assume that the metric on $M$ is selected to be a product near
$\partial M$. For $x\notin\partial M\times I$, let $R(x,\epsilon)$ be
the set of vectors in $T_x(M)$ of length less than $\epsilon$. If
$x\eq (y,t)\in\partial M\times I$, give $T_x(M)$ coordinates
$\omega_1,\ldots\,$,~$\omega_n$ so that
$\omega_1,\ldots\,$,~$\omega_{n-1}$ are in the $\partial M$-direction
(and hence exponentiate into $\partial M\times\set{t}$), and
$\omega_n$ is the coordinate in the $I$-direction. Then, define
$R(x,\epsilon)$ to be $\set{\omega\eq (\omega_1,\ldots, \omega_n)\in
T_x(M)\vbar\norm{\omega}<\epsilon\hbox{\ and\ }\omega_n\geq -t}$. For
small $\epsilon$ the exponential map $\Exp$ carries $R(x,\epsilon)$
diffeomorphically to an open neighborhood of $x$, even
when~$x\in\partial M$. For a properly imbedded submanifold $V$ of $M$
which meets $\partial M\times I$ in $I$-fibers, define
$N_\epsilon(V)\subset T(M)\vert_V$ by $N_\epsilon(V)\eq \cup_{x\in
V}R(x,\epsilon)$. When $V$ is compact, there exists a positive
$\epsilon$ such that for every $x\in V$, $\Exp$ carries each
$N_\epsilon(V)\cap T_x(M)$ diffeomorphically to a neighborhood of
$x\in M$.

For spaces of imbeddings, a ``$<\delta$'' subscript indicates the
imbeddings that are $\delta$-close to the inclusion, thus for example
$$\Imb_{<\delta}(V,M)\eq\set{j\in\Imb(V,M)\vbar
d(j(x),i_V(x))<\delta\hbox{\ for all $x\in V$}}\ .$$

\begin{lemma}{logarithm}{(Palais' lemma d)}
Assume that the metric on $M$ is a product near $\partial M$. Let $V$
be a compact submanifold of $M$ that meets $\partial M\times I$ in
$I$-fibers. For sufficiently small positive $\delta$, there exists a
continuous map $X\colon\Imb_{<\delta}(V,M)\to{\cal X}_{<1/2}(V,T(M))$
such that $\Exp(X(j))(x)\eq j(x)$ for all $x\in V$ and $j\in
\Imb_{<\delta}(V,M)$. Moreover, if $j(x)\eq i_V(x)$ then $X(j)(x)\eq
Z(x)$.
\marginwrite{logarithm}
\end{lemma}

\begin{proof}{} Choose $\epsilon<1/2$ small enough so that for all $x\in
V$, $\Exp$ gives a diffeomorphism from $N_\epsilon(V)\cap T_x(M)$ to a
neighborhood of $x$ in $M$. Choose $\delta$ small enough so that if
$j\in\Imb_{<\delta}(V,M)$ then $j(x)\in \Exp(N_\epsilon(V)\cap
T_x(M))$. For $j\in\Imb_{<\delta}(V,M)$ define $X(j)(x)$ to be the
unique vector in $N_\epsilon(V)\cap T_x(M)$ such that $\Exp(X(j)(x))$
equals~$j(x)$. We must verify that $X\in{\cal X}_{<1/2}(V,T(M))$.
Property~(1) holds because the metric is a product near
$\partial M$, so for $x\in\partial M$ and short vectors $\omega\in T_x(M)$,
$\Exp(\omega)\in\partial M$ if and only if $\omega$ is tangent
to~$\partial M$. Property~(2) and the final sentence of the lemma are
immediate.
\end{proof}

Before giving the main results of this section, we fix some notation
to simplify their statements. Suppose $V$ is a compact submanifold of
$M$, and $S\subseteq\partial M$ is a (possibly empty) closed subset
which is a neighborhood in $\partial M$ of $S\cap \partial V$. Note
that this implies that $S\cap\partial V$ is a union of components of
$V\cap\partial M$. In this situation, $\Imb(V,M\rel S)$ will stand for
the space of imbeddings that equal the inclusion on $V\cap S$ and
carry $V\cap (\partial M-S)$ into $\partial M-S$. As usual,
$\Imb^L(V,M\rel S)$ denotes the subspace consisting of all $j$ such
that $j(V)$ lies in the topological interior of~$L$.

The fundamental result is the analogue of theorem~B of \cite{P}. For
its proof we make one more definition. Define $F\colon{\cal
X}^L(T(M))\to \Maps^L(M,M)$ by $F(X)(x)\eq\Exp(X(x))$. We recall that
condition~(1) of the definition of ${\cal X}(T(M))$ and the fact that
the metric is a product near the boundary guarantee that
$F(X)(\partial M)\subset \partial M$. Since $\Exp$ is smooth, it
follows as in lemma~a of \cite{P} that $F$ is continuous.

\begin{theorem}{palaistheoremB}{} Let $V$ be a compact submanifold of
$M$, and let $S\subseteq\partial M$ be a closed neighborhood in
$\partial M$ of $S\cap \partial V$. Let $L$ be a neighborhood of $V$
in $M$. Then $\Imb^L(V,M\rel S)$ admits local $\Diff^L(M\rel S)$
cross-sections.
\marginwrite{palaistheoremB}
\end{theorem}

\begin{proof}{} By proposition~\ref{inclusion} it suffices to find a
local cross-section at the inclusion map $i_V$. Fix a compact
neighborhood $K$ of $V$ with $K\subseteq L$. Using
lemmas~\ref{logarithm} and~\ref{extension}, we obtain
$X\colon\Imb_{<\delta}(V,M)\to {\cal X}_{<1/2}(V,T(M))$ and $k\colon
{\cal X}_{<1/2}(V,T(M))\to {\cal X}^K(T(M))$. Let $J$ be a
neighborhood of $1_M$ in $\Maps^K(M,M)$ as in lemma~\ref{J}, and
define $U\eq (FkX)^{-1}(J)$. Then $\chi\eq FkX\colon U\to\Diff(M)$ is
the desired local $\Diff^L(M\rel S)$ cross-section at~$i_V$.
\end{proof}

From proposition~\ref{theoremA} we have immediate corollaries.

\begin{corollary}{palaiscoro2}{} Let $V$ be a compact submanifold of
$M$. Let $S\subseteq\partial M$ be a closed neighborhood in $\partial
M$ of $S\cap \partial V$, and $L$ a neighborhood of $V$ in $M$. Then
the restriction $\Diff^L(M\rel S)\to \Imb^L(V,M\rel S)$ is locally
trivial.
\marginwrite{palaiscoro2}
\end{corollary}

\begin{corollary}{palaiscoro3}{} Let $V$ and $W$ be compact
submanifolds of $M$, with $W\subseteq V$. Let $S\subseteq\partial M$ a
closed neighborhood in $\partial M$ of $S\cap \partial V$, and $L$ a
neighborhood of $V$ in $M$. Then the restriction $\Imb^L(V,M\rel S)
\to\Imb^L(W,M\rel S)$ is locally trivial.
\marginwrite{palaiscoro3}
\end{corollary}

\section{The vertical and aligned exponentials}
\label{exponent}\marginwrite{exponent}

Let $p\colon E\to B$ be a locally trivial smooth map of manifolds,
with compact fiber, and let $\pi\colon T(E)\to E$ denote the tangent
bundle of $E$. At each point $x\in E$, let $V_x(E)$ denote the {\it
vertical subspace} of $T_x(E)$ consisting of vectors tangent to the
fiber of~$p$. When $E$ has a Riemannian metric, the orthogonal
complement $H_x(E)$ of $V_x(E)$ in $T_x(E)$ is called the {\it
horizontal subspace.} We usually abbreviate $V_x(E)$ and $H_x(E)$ to
$V_x$ and $H_x$, and call their elements {\it vertical} and {\it
horizontal} respectively. Clearly $V_x$ is the kernel of $p_*\colon
T_x(E)\to T_{p(x)}(B)$, while $p_*\vert_{H_x}\colon H_x\to
T_{p(x)}(B)$ is an isomorphism. Each vector $\omega\in T_x(E)$ has an
orthogonal decomposition $\omega\eq \omega_v+\omega_h$.

Define the {\it horizontal boundary} $\partial_hE$ to be $\cup_{x\in
B}\partial(p^{-1}(x))$, and the {\it vertical boundary} $\partial_vE$
to be $p^{-1}(\partial B)$.

A path $\alpha$ in $E$ is called {\it horizontal} if $\alpha'(t)\in
H_{\alpha(t)}$ for all $t$ in the domain of $\alpha$. Let
$\gamma\colon [a,b]\to B$ be a path such that $\gamma'(t)$ never
vanishes, and let $x\in E$ with $p(x)\eq \gamma(a)$. A horizontal path
$\widetilde{\gamma}\colon [a,b]\to E$ such that
$\widetilde{\gamma}(a)\eq x$ and $p\widetilde{\gamma}\eq \gamma$ is
called a {\it horizontal lift} of $\gamma$ starting at~$x$.

To ensure that horizontal lifts exist, we will need a special metric
on $E$. Using the local product structure, at each point $x$ in
$\partial_hE$ select a vector field defined on a neighborhood of
$x$ that
\begin{enumerate}
\item[{\rm(a)}] points into the fiber at points of $\partial_hE$, and
\item[{\rm(b)}] is tangent to the fibers wherever it is defined.
\end{enumerate}

\noindent Note that by (b), the vector field must be tangent to
$\partial_vE$ at points in $\partial_vE$.  Since scalar multiples and
linear combinations of vectors satisfying these two conditions also
satisfy them, we may piece these local fields together using a
partition of unity to construct a vector field, nonvanishing on a
neighborhood of $\partial_hE$, that satisfies (a) and~(b). Using the
integral curves associated to this vector field we obtain a smooth
collar neighborhood $\partial_hE\times [0,2]$ of $\partial_hE$ such
that each $[0,2]$-fiber lies in a fiber of $p$. On $\partial_hE\times
[0,2)$, fix a Riemannian metric that is the product of a metric on
$\partial_hE$ and the usual metric on $[0,2)$. Form a metric on $E$
from this metric and any metric on all of $E$ using a partition of
unity subordinate to the open cover $\set{\partial_hE\times[0,2),
E-\partial_hE\times I}$. Such a metric is said to be a {\it product
near $\partial_hE$ such that the $I$-fibers of $\partial_hE\times I$
are vertical.}  It has the following properties for $0\leq t\leq 1$:
\begin{enumerate}
\item[{\rm(i)}] If $x\in E-\partial_hE\times I$, then $d(x,\partial_hE)>1$.
\item[{\rm(ii)}] If $x\eq(y,t)\in\partial_hE\times\set{t}$, then
$d(x,\partial_hE)\eq t$.
\item[{\rm(iii)}] For $x\in\partial_hE\times\set{t}$, the horizontal
subspace $H_x$ is tangent to~$\partial_hE\times\set{t}$.
\end{enumerate}

\noindent To see property (iii), start with the fact that $H_x$ is
perpendicular to the fiber $p^{-1}(p(x))$. Since the $I$-fiber of
$\partial_hE\times I$ that contains $x$ lies in $p^{-1}(p(x))$, $H_x$
is orthogonal to that $I$-fiber as well. Since
$\partial_hE\times\set{t}$ meets the $I$-fiber orthogonally, with
codimension~1, $H_x$ is tangent to $\partial_hE\times\set{t}$.

Property (iii) implies that a horizontal lift starting in some
$\partial_hE\times \set{t}$ will continue in $\partial_hE\times
\set{t}$. Using the compactness of the fiber, the existence of horizontal
lifts will then be guaranteed.

\begin{lemma}{hlift}{}
Assume that the metric on $E$ is a product near $\partial_hE$ such
that the $I$-fibers of $\partial_hE\times I$ are vertical. Let
$\gamma\colon [a,b]\to B$ be a path such that $\gamma'(t)$ never
vanishes, and let $x\in E$ with $p(x)\eq \gamma(a)$. Then there exists
a unique horizontal lift of $\gamma$ starting at~$x$.
\marginwrite{hlift}
\end{lemma}

\begin{proof}{} For any horizontal lift $\widetilde{\gamma}(t)$, each
$\widetilde{\gamma}'(t)$ is uniquely determined, so the lift through a
given point in $E$ is unique if it exists. For each $\gamma(t)$, let
$F_{\gamma(t)}$ be the fiber over $\gamma(t)$. From the local theory
of ordinary differential equations, each point in $F_{\gamma(t)}$ that
does not lie in $\partial_hE$ has a neighborhood in
$p^{-1}(\gamma([a,b])$ in which $\gamma$ has horizontal lifts. Since
the metric is a product near $\partial_hE$ such that the $I$-fibers
are vertical, the same is true for each point in $\partial_hE$.  Since
the fiber is compact, for each $t$ there exists an $\epsilon$ such
that for every $x\in F_{\gamma(t)}$, the horizontal lift of $\gamma$
through $x$ exists for $s\in(t-\epsilon,t+\epsilon)$, and the result
follows using compactness of the interval~$[a,b]$.
\end{proof}

\noindent For the remainder of this section, assume that the
metric on $E$ is a product near $\partial_hE$ such that the $I$-fibers
of $\partial_hE\times I$ are vertical. Each fiber $F$ of $E$ inherits
a Riemannian metric from that of $E$, and has an exponential map
$\Exp_F$ which (where defined) carries vectors tangent to $F$ to
points of $F$. Note that the path $\Exp_F(t\omega)$ is not generally a
geodesic in $E$. The vertical exponential $\Exp_v$ is defined by
$\Exp_v(\omega)\eq\Exp_F(\omega)$, where $\omega$ is a vertical vector
and $F$ is the fiber containing~$\pi(\omega)$.

Before defining the aligned exponential map $\Exp_a$,
we will motivate its definition. A vector field $X\colon E\to T(E)$ is
called {\it aligned} if $p(x)\eq p(y)$ implies that $p_*(X(x))\eq
p_*(X(y))$. This happens precisely when there exist a vector field
$X_B$ on $B$ and a vertical vector field $X_V$ on $E$ so that for all
$x\in E$,
$$X(x)\eq (p_*\vert_{H_x})^{-1}(X_B(p(x)))+X_V(x)\ .$$

\noindent In particular, any vertical vector field is aligned. When
$X$ is aligned, the projected vector field $p_*X$ is well-defined.
The key property of $\Exp_a$ is that if $X$ is an aligned vector field
on $E$, and $\Exp_a(X(x))$ is defined for all $x$, then the map of $E$
defined by sending $x$ to $\Exp_a(X(x))$ will be fiber-preserving.

$\Exp_a$ is defined as follows. Consider a tangent vector $\beta$ in
$B$ such that $\Exp(\beta)$ is defined.  A geodesic segment
$\gamma_\beta$ starting at $\pi(\beta)$ is defined by
$\gamma_\beta(t)\eq \Exp(t\beta)$, $0\leq t\leq 1$.  Define
$\Exp_a(\omega)$ to be the endpoint of the unique horizontal lift of
$\gamma_{p_*(\omega)}$ starting at $\Exp_v(\omega_v)$. Note that
$\Exp_a(\omega)$ exists if and only if both $\Exp_v(\omega_v)$ and
$\Exp(p_*(\omega))$ exist.  Clearly, when $\Exp_a(\omega)$ is defined,
it lies in the fiber containing the endpoint of a lift of
$\gamma_{p_*(\omega)}$, and therefore $p(\Exp_a(\omega))\eq
\Exp(p_*(\omega))$. This immediately implies that if $X$ is an aligned
vector field on $E$ such that $\Exp_a(X(x))$ is defined for all $x\in
E$, then the map defined by sending $x$ to $\Exp_a(X(x))$ takes fibers
to fibers, and in particular if $X$ is vertical, it takes each fiber
to itself.

In section~\ref{restrict} we will need a further refinement of the
metric on $E$, namely that it also be a product near $\partial_vE$. To
achieve this, we proceed as follows. If $\partial_vE$ is empty, there
is nothing needed, and if $\partial_hE$ is empty, then we simply
choose a metric which is a product near the boundary as in
section~\ref{metrics}. Assuming that both are nonempty, put $Y\eq
\partial(\partial_vE)\eq \partial(\partial_hE)\eq
\partial_vE\cap \partial_hE$. Let $R_h$ be the complete metric on
$\partial_hE$ that was was used to construct the metric $R$ on $E$
that is a product on a collar $\partial_hE\times[0,1]_1$, where the
subscript will distinguish this interval from another to be selected
later. Since the choice of $R_h$ was arbitrary, we may assume that
$R_h$ was a product near $Y$. That is, after reparametrizing, we may
assume that there is a collar $Y\times[0,2]_2$ of $Y$ in $\partial_hE$
such that $R_h$ is a product on all of $Y\times[0,2]_2$. Now
$Y\times[0,2]_1$ is a collar of $Y$ in $\partial_vE$, and
$Y\times[0,2]_1\times[0,2]_2$ is a partial collar of $\partial_vE$
defined on the subset $Y\times[0,2]_1$. Extend this to a collar
$\partial_vE\times[0,2]_2$.  Let $R_v$ be the product of the
restricted metric $R\vert_{\partial_vE}$ and the standard metric on
$[0,2]_2$.  Since $R_h$ was a product near $Y$, we have $R_v\eq R$ on
$Y\times[0,1]_1\times [0,2]_2$. Now form a new metric on $E$ by
piecing together $R_v$ and $R$ using a partition of unity subordinate
to the open cover $\set{\partial_vE\times[0,2)_2,
E-\partial_vE\times[0,1]_2}$. On points of
$Y\times[0,1]_1\times[0,2)_2$ the new metric is just a linear
combination of the form $tR+(1-t)R$, so agrees with $R$. The resulting
metric is both a product near $\partial_hE$ and a product
near~$\partial_vE$.

\section{Projection of fiber-preserving diffeomorphisms}
\label{project}\marginwrite{project}

Throughout this section and the next, it is understood that $p\colon
E\to B$ is a locally trivial smooth map as in section~\ref{exponent},
such that the metric on $B$ is a product near $\partial B$, and the
metric on $E$ is a product near $\partial_hE$ such that the $I$-fibers
of $\partial_hE\times I$ are vertical. When $W$ is a vertical
submanifold of $E$, it is automatic that $W$ meets the collar
$\partial_hE\times I$ in $I$-fibers, and we by rechoosing the metric
on $B$ we may assume that $p(W)$ meets the collar $\partial B\times I$
in $I$-fibers. Define $\partial_hW\eq W\cap\partial_hE$ and
$\partial_vW\eq W\cap\partial_vE$.

Define ${\cal A}(W,T(E))$ to be the sections $X$ from $W$ to
$T(E)\vert_W$ such that
\begin{enumerate}
\item[(1)] $X$ is aligned, that is, if $p(w_1)\eq p(w_2)$ then
$p_*(X(w_1))\eq p_*(X(w_2))$,
\item[(2)] if $x\in \partial_hW$, then $X(x)$ is tangent to
$\partial_hE$, and if $x\in \partial_vW$, then $X(x)$ is tangent to
$\partial_vE$, and
\item[(3)] $\Exp_a(X(x))$ is defined for all $x\in W$.
\end{enumerate}

\noindent
When $W\eq E$, the vector fields satisfying~(1), (2), and~(3) are
denoted by ${\cal A}(T(E))$. The embellishments ${\cal A}^L(T(E))$ and
${\cal A}_{<1/2}(W,T(E))$ have the same meanings as in
section~\ref{palais}. The elements of ${\cal A}(W,T(E))$ such that
$p_*X(x)\eq Z(p(x))$ for all $x\in W$ are denoted by ${\cal
V}(W,T(E))$, and similarly for ${\cal V}(T(E))$.

Define $F_a\colon{\cal A}^L(T(E))\to \Maps^L(E,E)$ by
$F_a(X)(x)\eq\Exp_a(X(x))$. Since $\Exp_a$ is smooth, it follows as in
lemma~a of \cite{P} that $F_a$ is continuous.

\begin{theorem}{theorem1}{} Let $K$ be a compact subset of $B$. Let
$S$ be a subset of $\partial B$, and let $T\eq p^{-1}(S)$.  Then
$\Diff^K(B\rel S)$ admits local $\Diff^{p^{-1}(K)}_f(E\rel T)$
cross-sections.
\marginwrite{theorem1}
\end{theorem}

\begin{proof}{} Choose a compact subset $L$ of $B$ such that $K\subseteq
\hbox{\it int}(L)$. By lemma~\ref{logarithm}, there exist $\delta>0$
and a continuous map $X_1\colon \Imb_{<\delta}(L,B)\to {\cal
X}_{<1/2}(L,T(B))$ such that $\Exp(X_1(j)(x))\eq j(x)$ for all $x\in
L$ and all $j\in\Imb_{<\delta}(L,B)$. Moreover, if $j(x)\eq x$, then
$X_1(j)(x)\eq Z(x)$.

Let $\rho\colon \Diff^K(B)\to \Imb(L,B)$ be restriction. Put $U_0\eq
\rho^{-1}(\Imb_{<\delta}(L,B))$, a neighborhood of $1_B$, and define
$X_0\colon U_0\to {\cal X}(T(B))$ by
$$X_0(f)(x)\eq\cases{X_1(f\vert_L)(x)&if $x\in L$\cr
                      Z(x)&if $x\notin K$\ .\cr}$$

\noindent This makes sense since if $x\in L-K$, then $f\vert_L(x)\eq
x$ so $X_1(f\vert_L)(x)\eq Z(x)$. We have $X_0(1_B)\eq Z$ and
$\Exp(X_0(f)(x))\eq f(x)$ for every $f\in U_0$ and $x\in B$.

Let $h\in \Diff^K(B)$. For every $g\in U_0h$, $\Exp(X_0(gh^{-1}(x)))\eq
gh^{-1}(x)$. Define
$$\widetilde{\chi}(g)(x)=
\big(p_*\vert_{H_x}\big)^{-1}(X_0(gh^{-1})(p(x)))\ ,$$

\noindent so that $\widetilde{\chi}(g)$ is an aligned section of
$T(E)$. We have that $\Exp_a(\widetilde{\chi}(g)(x))$ exists since
$\Exp(p_*\widetilde{\chi}(g)(x))\eq \Exp(X_0(gh^{-1})(p(x)))\eq
gh^{-1}(p(x))$ exists and $\Exp_v(\widetilde{\chi}(g)(x))\eq x$
exists. The other conditions are easily checked to verify that
$\widetilde{\chi}(g)\in{\cal A}^{p^{-1}(K)}(T(E))$.

Let $J$ be a neighborhood of $1_M\in\Maps^{p^{-1}(K)}(E,E)$ given by
lemma~\ref{J}, and put $U\eq\widetilde{\chi}^{-1}F_a^{-1}(J)$. Define
$\chi\colon U\to\Diff_f(E)$ by $\chi(g)\eq F_a\widetilde{\chi}(g)$.
The local cross-section condition holds, since given $b\in B$ we may
choose $x$ with $p(x)\eq h(b)$ and calculate that
$$\eqalign{\ol{\chi(g)}h(b)&= p(\chi(g)(x))\cr
&=p\big(\Exp_a\big(\widetilde{\chi}(g)(x)\big)\big)\cr
&=\Exp\big(X_0(gh^{-1})(h(b))\big)\cr
&=gh^{-1}(h(b))\cr
}$$

\noindent If $g,h\in\Diff^K(B\rel S)$, then
$X_0(gh^{-1})(x)\eq Z(x)$ for all $x\in\partial S$. It follows that
$\widetilde{\chi}(g)(x)\eq Z(x)$ for all $x\in T$, so
$\chi(g)\in\Diff^{p^{-1}(K)}_f(E\rel T)$.
\end{proof}

From proposition~\ref{theoremA}, we have immediately

\begin{theorem}{project diffs}{} Let $K$ be a compact subset of $B$.
Let $S\subseteq\partial B$ and let $T\eq p^{-1}(T)$. Then
$\Diff^{p^{-1}(K)}_f(E\rel T)\to \Diff^K(B\rel S)$ is locally trivial.
\marginwrite{project diffs}
\end{theorem}

The homotopy extension property of the fibration in
theorem~\ref{project diffs} yields immediately the following
corollary. As indicated in the introduction, each of the fibration
theorems we prove in this paper has a corresponding corollary
involving parameterized lifting or extension, but since the statements
are all analogous we give only the following one as a prototype.

\begin{corollary}{isotopy lifting}{(Parameterized Isotopy Extension
Theorem)} Let $K$ be a compact subset of $B$, let $S\subseteq\partial
B$, and let $T\eq p^{-1}(S)$. Suppose that for each $t$ in a
path-connected parameter space $P$ there is an isotopy $g_{t,s}$,
which is the identity on $S$ and outside of $K$, such that $g_{t,0}$
lifts to a diffeomorphism $G_{t,0}$ of $E$ which is the identity on
$T$. Assume that sending $(t,s)\to g_{t,s}$ defines a continuous
function from $P\times [0,1]$ to $\Diff(B\rel S)$ and sending $t$ to
$G_{t,0}$ defines a continuous function from $P$ to $\Diff(E\rel
T)$. Then the family $G_{t,0}$ extends to a continuous family on
$P\times I$ such that for each $(t,s)$, $G_{t,s}$ is a
fiber-preserving diffeomorphism inducing $g_{t,s}$ on~$B$.
\marginwrite{isotopy lifting}
\end{corollary}

\section{Restriction of fiber-preserving diffeomorphisms}
\label{restrict}
\marginwrite{restrict}

In this section we present the analogues of the main results of
\cite{P} in the fibered case. As in section~\ref{project}, we
assume that the metric on $B$ is a product near $\partial B$, and the
metric on $E$ is a product near $\partial_hE$ such that the $I$-fibers
of $\partial_hE\times I$ are vertical. We further assume that the
metric on $E$ is a product near $\partial_vE$; this is needed only in
the proof of lemma~\ref{lemmaC}.

It is first necessary to adapt lemmas~\ref{extension} and~\ref{logarithm}.

\begin{lemma}{lemmaC}{}
Let $W$ be a compact vertical submanifold of $E$. Let $T$ be a closed
fibered neighborhood in $\partial_vE$ of $T\cap\partial_vW$, and let
$L\subseteq E$ be a neighborhood of $W$. For sufficiently small
$\delta$ there exists a continuous map $k\colon{\cal
A}_{<\delta}(W,T(E))\to {\cal A}^L(T(E))$ such that $k(X)(x)\eq X(x)$
for all $x\in W$ and $X\in {\cal A}_{<\delta}(W,T(E))$. If $X(x)\eq
Z(x)$ for all $x\in T\cap\partial_vW$, then $k(X)(x)\eq Z(x)$ for all
$x\in T$. Furthermore, $k({\cal V}_{<\delta}(W,T(E)))\subset {\cal
V}^L(T(E))$.
\marginwrite{lemmaC}
\end{lemma}

\begin{proof}{} Assume first that $W$ has positive codimension. Since
the fiber of $p$ is compact, we may assume that $p(L)$ is a
neighborhood of $p(W)$ with $p^{-1}(p(L))\eq L$. Since $W$ is compact
we may choose $\delta<1/2$ such that if $X\in {\cal
A}_{<\delta}(W,T(E))$ then $\norm{p_*X(x)}<1/2$ for all $x\in p(W)$.
Let $k_B\colon {\cal X}_{<1/2}(p(W),T(B)) \to{\cal X}^{p(L)}(T(B))$ be
given by lemma~\ref{extension}, using $p(T)$ as the neighborhood $S$
in lemma~\ref{extension}.  The vectors
$(p_*\vert_{H_x})^{-1}(k_B(p_*X)(p(x)))$ will give the horizontal part
of our extension $k$, but to obtain sufficient control of the vertical
part we will need to adapt the proof of lemma~\ref{extension}
using~$\Exp_a$.

Let $\nu_\epsilon(W)$ be the $\epsilon$-normal bundle of $W$. Note
that its fibers are horizontal, since the tangent space of $W$
includes the bundle of vertical vectors of $T(E)$ at points of $W$.
For sufficiently small~$\epsilon$, $j_a\colon \nu_\epsilon(W)\to E$
can be defined by $j_a(\omega)\eq\Exp_a(\omega)$ and carries
$\nu_\epsilon(W)$ diffeomorphically to a neighborhood of $W$ in $E$.
We choose $\epsilon$ small enough so that this neighborhood is
contained in $L$, and so that the image of the fibers at points of
$T\cap\partial_vW$ lies in $T$ and the image of the fibers at points
of $(\partial_vE-T)\cap\partial_vW$ lies in~$\partial_vE-T$.

If $x\in\partial_hE\times \set{t}$, $j_a$ carries the normal fiber at
$x$ into $\partial_hE\times\set{t}$.  Since $W$ is compact, we may
choose $\epsilon$ small enough so that $j_a(\omega)\in
\partial_hE\times I$ only when $\pi(\omega)\in\partial_hE\times I$.

Suppose $v\in T_x(E)$ and that $\Exp_a(v)$ is defined. For all $u\in
T_x(E)$ define $P_a(u,v)$ to be the vector that results from parallel
translation of $u$ along the path that sends $t$ to~$\Exp_a(tv)$, $0\leq
t\leq 1$. Let $\alpha\colon E\to [0,1]$ be a smooth function which is
identically~1 on $W$ and identically~0 on
$E-j_a(\nu_{\epsilon/2}(W))$. Define $k_E\colon {\cal
A}_{<\delta}(W,T(E))\to {\cal V}^L(T(E))$ by
$$k_E(X)(x)\eq
\cases{
\alpha(x)P_a(X(\pi(j_a^{-1}(x))),j_a^{-1}(x))_v&for
$x\in j_a(\nu_{\epsilon}(W))$\cr
Z(x)&for $x\in E-j_a(\nu_{\epsilon/2}(W))$\ .\cr}$$

\noindent Note that if $X$ is vertical, then $k(X)(x)\eq X(x)$ for all
$x\in W$. For later use we make two observations about $k_E$.  First,
if $X(x)\eq 0$ for all $x\in T\cap\partial_vW$, then $k_E(X)(x)\eq 0$
for all $x\in T$. This is because $j_a^{-1}(T)$ consists exactly of
the vectors normal to $W$ at the points of $T\cap\partial_vW$. Second,
if $x\eq(y,t)\in\partial_hE\times I$, and $x\in j_a(\nu_\epsilon(W))$,
then $\pi(j_a^{-1}(x))$ is of the form $(y',t)$, and either
$k_E(X)(x)\eq Z(x)$ or the component of $k_E(X)(x)$ in the
$I$-direction is of the form $\beta\omega_I$, where $0<\beta\leq1$ and
$\omega_I$ is the component of $X(\pi(j_a^{-1}(x)))$ in the
$I$-direction. This follows because the metric is a product on
$\partial_hE\times I$, so parallel translation preserves the component
in the $I$-direction. Consequently, since
$\Exp_v(X(\pi(j_a^{-1}(x))))$ is defined, so is $\Exp_v(k(X)(x))$.

For $X\in{\cal A}_{<\delta}(W,T(E))$, define $X_v(x)=X(x)_v$, and put
$$k(X)(x)=(p_*\vert_{H_x})^{-1}(k_B(p_*X)(p(x)))\;+\;k_E(X_v)(x)\ .$$

\noindent
We need to check that $k(X)$ lies in ${\cal A}^L(T(E))$. From its
definition, $k(X)$ is aligned and vanishes outside of $L$. Let
$x\in\partial_hE$ and suppose that $k(X)(x)\neq Z(x)$. Then
$\pi(j_a^{-1}(x))\in\partial_hE$ and since $X\in{\cal A}(W,T(E))$,
$X(\pi(j_a^{-1}(x)))$ is tangent to $\partial_hE$. Since the metric is
a product near $\partial_hE$, parallel translation preserves vectors
tangent to $\partial_hE$, and it follows that $k(X)(x)$ is tangent to
$\partial_hE$. Suppose that $x\in\partial_vE$ and $k(X)(x)\neq Z(x)$.
Since $k_B(p_*X)(p(x))$ is tangent to $\partial B$,
$\big(p_*\vert_{H_x}\big)^{-1}(k_B(p_*X)(p(x)))$ is tangent to
$\partial_vE$. The fact that the metric on $E$ is a product near
$\partial_vE$ implies that $\pi(j_a^{-1}(x))\in\partial_vW$, and
moreover, since $X(\pi(j_a^{-1}(x)))$ is tangent to $\partial_vE$,
that $P_a(X(\pi(j_a^{-1}(x))),j_a^{-1}(x))$ is also tangent to
$\partial_vE$. We conclude that $k(X)(x)$ is tangent to $\partial_vE$.
The fact that $\Exp(k_B(p_*X)(p(x)))$ was defined, together with the
second observation after the definition of~$k_E$, implies that
$\Exp(k(X)(x))$ is always defined.

Suppose that $X(x)\eq Z(x)$ for all $x\in T\cap\partial_vW$. Then
$p_*(X)(p(x))\eq Z(p(x))$ for all $x\in p(T)\cap\partial(p(W))$, so
$k_B(p_*X)(p(x))\eq Z(p(x))$ for all $x\in p(T)$. The first
observation after the definition of $k_E$ shows that $k_E(X)(x)\eq
Z(x)$ for all $x\in T$. Therefore $k(X)(x)\eq Z(x)$ for all~$x\in T$.

For the last statement, if $X\in{\cal V}(W,T(E))$, then
$p_*(X)(p(x))\eq Z(p(x))$ for all $x\in \partial(p(W))$, so
$k_B(p_*X)(p(x))\eq Z(p(x))$ for all $x\in\partial_vE$. Therefore
$k(X)\in{\cal V}(T(E))$.

The case when $W$ has codimension zero is similar. As in the proof of
lemma~\ref{extension}, use the subset $\nu_\epsilon^+\Fr(W)$
consisting of the vectors in the normal bundle of the frontier of $W$
whose aligned exponential lies in $\overline{E-W}$, and define
$$k_E(X)(x)\eq
\cases{X(x)&for $x\in W$\cr
\alpha(x)P_a(X(\pi(j_a^{-1}(x))),j_a^{-1}(x))_v&for
$x\in j_a(\nu^+_{\epsilon}(\Fr(W)))$\cr
Z(x)&for $x\in E-j_a(\nu^+_{\epsilon/2}(\Fr(W)))$\ .\cr}$$
\end{proof}

For the next lemma we will adapt the neighborhood $N_\epsilon(V)$ used
in the proof of lemma~\ref{logarithm} into the fibered context. For
$x\in E$, let $R_B(p(x),\epsilon)$ be the subset of $T_{p(x)}(B)$ as
defined before the statement of lemma~\ref{logarithm}. Denote
$p^{-1}(p(x))$ by $F$, and let $R_F(x,\epsilon)$ be the subset of
$T_x(F)$ defined before the statement of lemma~\ref{logarithm}. Regard
$T_x(F)$ as the vertical subset $V_x(E)$ of $T_x(E)$, and observe that
for sufficiently small $\epsilon$, the aligned exponential $\Exp_a$
carries $R_F(x,\epsilon)$ to a neighborhood of $x$ in $R$, since on
$T_x(F)$, $\Exp_a$ agrees with the exponential map of~$F$. Now define
$S(x,\epsilon)\eq R_F(x,\epsilon)\times
\Big(p_*\vert_{H_x}\Big)^{-1}(R_B(p(x),\epsilon))\subset V_x(E)\times
H_x(E)\eq T_x(E)$. For a vertical submanifold $W\subseteq E$, define
$N_\epsilon(W)\eq \cup_{x\in W}S(x,\epsilon)$. Provided that $W$ is
compact, we may choose a positive $\epsilon$ such that for each $x\in
W$, $\Exp_a\colon N_\epsilon(W)\cap T_x(E)\to E$ is a diffeomorphism
onto a neighborhood of $x$ in~$E$.

\begin{lemma}{lemmaD}{}
Let $W$ be a compact vertical submanifold of $E$.  For sufficiently
small $\delta$, there exists a continuous map $X\colon
(\Imb_f)_{<\delta}(W,E)\to{\cal A}(W,T(E))$ such that
$\Exp_a(X(j)(x))\eq j(x)$ for all $x\in W$ and $j\in
(\Imb_f)_{<\delta}(W,E)$, and moreover if $j(x)\eq i_W(x)$ then
$X(j)(x)\eq Z(x)$. Furthermore, $X((\Imb_v)_{<\delta}(W,E))\subset
{\cal V}(W,T(E))$.
\marginwrite{lemmaD}
\end{lemma}

\begin{proof}{} Let $N_\epsilon(W)$ be as defined above,
with $\epsilon$ small enough to ensure the local diffeomorphism
condition.  Choose $\delta$ small enough so that for
every $x\in W$ and every $j\in (\Imb_f)_{<\delta}(W,E)$,
$j(x)\in\Exp_a(N_\epsilon(W)\cap T_x(E))$. Define $X(j)(x)$ to be the
unique vector in $N_\epsilon(W)\cap T_x(M)$ such that $\Exp_a(X(j)(x))$
equals~$j(x)$.
\end{proof}

In the statements of our remaining results, the notation $\Imb(W,E\rel
T)$ is as defined before the statement of
theorem~\ref{palaistheoremB}.

\begin{theorem}{theorem2}{}
Let $W$ be a compact vertical submanifold of $E$. Let $T$ be a closed
fibered neighborhood in $\partial_vE$ of $T\cap \partial_vW$, and let
$L$ be a neighborhood of $W$. Then
\begin{enumerate}
\item[{\rm (i)}] $\Imb_f^L(W,E\rel T)$ admits local $\Diff_f^L(E\rel
T)$ cross-sections, and
\item[{\rm (ii)}] $\Imb_v(W,E\rel T)$ admits local $\Diff_v^L(E\rel
T)$ cross-sections.
\end{enumerate}
\marginwrite{theorem2}
\end{theorem}

\begin{proof}{} By proposition~\ref{inclusion}, it suffices to find
local cross-sections at the inclusion $i_W$. Choose a compact
neighborhood $K$ of $W$ with $K\subseteq L$. Let $k\colon{\cal
A}_{<\delta}(W,T(E))\to{\cal X}^K(T(E))$ be obtained using
lemma~\ref{lemmaC}. Using lemma~\ref{lemmaD}, choose $\delta_1>0$ and
$X\colon (\Imb_f)_{<\delta_1}(W,E)\to {\cal A}(W,T(E))$. If $j\in
\Imb_f^L(W,E\rel T)$, then $X(j)(x)\eq 0$ for all $x\in\partial_vE$.
Since $X$ is continuous and $X(i_W)\eq Z$, we may choose a
neighborhood $U$ of $i_W$ in $\Imb_f^L(W,E\rel T)$ so that
$X(U)\subset {\cal A}_{<\delta}(W,T(E))$. For $j\in U$ define
$\chi(j)\eq F_akX(j)$. From lemma~\ref{J}, we may make $U$ small
enough to ensure that $\chi(j)$ is a diffeomorphism, and then $\chi$
is the cross-section that proves~(a). For~(b), suppose that $j\in
U\cap\Imb_v^L(W,E\rel T)$. Since $k({\cal
V}_{<\delta}(W,T(E)))\subseteq{\cal V}^K(T(E))$, $\chi(j)$ lies in
$\Diff_v(E\rel T)$, so the restriction of $\chi$ to
$U\cap\Imb_v^L(W,E\rel T)$ is the necessary cross-section.
\end{proof}

Using proposition~\ref{theoremA}, we have the following immediate
corollaries.

\begin{corollary}{corollary2}{}
Let $W$ be a compact vertical submanifold of $E$. Let $T$ be a closed
fibered neighborhood in $\partial_vE$ of $T\cap \partial_vW$, and $L$
a neighborhood of $W$. Then the following restrictions are locally
trivial:
\begin{enumerate}
\item[{\rm(i)}] $\Diff_f^L(E\rel T)\to\Imb_f^L(W,E\rel T)$, and
\item[{\rm(ii)}] $\Diff_v^L(E\rel T)\to \Imb_v(W,E\rel T)$.
\end{enumerate}
\marginwrite{corollary2}
\end{corollary}

\begin{corollary}{corollary3}{} Let $V$ and $W$ be vertical
submanifolds of $E$. Let $T$ be a closed fibered neighborhood in
$\partial_vE$ of $T\cap \partial_vV$, and let $L$ a neighborhood of
$V$. Then the following restrictions are locally trivial:
\begin{enumerate}
\item[{\rm(i)}] $\Imb_f^L(V,E\rel T)\to\Imb_f^L(W,E\rel T)$.
\item[{\rm(ii)}] $\Imb_v(V,E\rel T)\to \Imb_v(W,E\rel T)$.
\end{enumerate}
\marginwrite{corollary3}
\end{corollary}

The final result of this section includes some of our previous
results.

\begin{theorem}{square}{} Let $W$ be a compact vertical submanifold of
$E$. Let $K$ be a compact neighborhood of $p(W)$ in $B$. Let $T$ be a
closed fibered neighborhood in $\partial_v\Sigma$ of $T\cap
\partial_vW$, and put $S\eq p(T)$. Then all four maps in the following
square are locally trivial:
$$\vbox{\halign{\hfil#\hfil\quad&#&\quad\hfil#\hfil\cr
$\Diff_f^{p^{-1}(K)}(E\rel T)$&$\longrightarrow$&%
$\Imb_f^{p^{-1}(K)}(W,E\rel T)$\cr
\noalign{\smallskip}
$\mapdown{}$&&$\mapdown{}$\cr
\noalign{\smallskip}
$\Diff^K(B\rel S)$&$\longrightarrow$&%
$\Imb^K(p(W),B\rel S)$\rlap{\ .}\cr}}$$
\marginwrite{square}
\end{theorem}

\begin{proof}{} The top arrow is corollary~\ref{corollary2}(i), the
left vertical arrow is theorem~\ref{theorem1}, and the bottom arrow is
corollary~\ref{palaiscoro2}.  For the right vertical arrow, we will
first show that $\Imb^K(p(W),B \rel\allowbreak S)$ admits local
$\Diff_f^{p^{-1}(K)}(E\allowbreak \rel\allowbreak T)$ cross-sections.
Let $i\in \Imb^K(p(W),B\rel S)$.  Using theorems~\ref{theorem2}
and~\ref{theorem1}, choose local cross-sections $\chi_1\colon U\to
\Diff^K(B\rel S)$ at $i$ and $\chi_2\colon
V\to\Diff_f^{p^{-1}(K)}\allowbreak(E\rel T)$ at $\chi_1(i)$.  Let
$U_1\eq\chi_1^{-1}(V)$, then for $j\in U_1$ we have
$$\overline{\chi_2\chi_1(j)}i=\overline{\chi_2(\chi_1(j))}i=
\chi_1(j)i=j\ .$$

\noindent Since the right vertical arrow is
$\Diff_f^{p^{-1}(K)}(E\rel T)$-equivariant, proposition \ref{theoremA}
implies it is locally trivial.
\end{proof}

\section{Palais' theorem for orbifolds}
\marginwrite{orbifold}
\label{orbifold}

In this section, we prove the main results from \cite{P} in the
context of orbifolds. Let ${\cal O}$ be a connected smooth orbifold
whose universal covering $\widetilde{\cal O}$ is a manifold. Denote by
$\tau\colon\widetilde{\cal O}\to{\cal O}$ the orbifold universal
covering, and let $H$ be its group of covering transformations. Since
${\cal O}$ is smooth, the elements of $H$ are diffeomorphisms.

Let $\Maps_H(\widetilde{\cal O},\widetilde{\cal O})$ be the space of
weakly $H$-equivariant maps, that is, the maps $f\colon
\widetilde{\cal O}\to \widetilde{\cal O}$ such that for some
automorphism $\alpha$ of $H$, $f(h(x))=\alpha(h)(f(x))$ for all $x\in
\widetilde{\cal O}$ and $h\in H$. Let $\Diff_H(\widetilde{\cal O})$ be
the weakly $H$-equivariant diffeomorphisms of $\widetilde{\cal O}$.
It is the normalizer of $H$ in $\Diff(\widetilde{\cal O})$.

An {\it orbifold diffeomorphism} of ${\cal O}$ is by definition an
orbifold homeomorphism of ${\cal O}$ whose lifts to $\widetilde{\cal
O}$ are diffeomorphisms. Thus the group $\Diff({\cal O})$ of orbifold
diffeomorphisms of ${\cal O}$ is the quotient of the group
$\Diff_H({\widetilde{\cal O}})$ by the normal subgroup~$H$.

Throughout this section, we let ${\cal W}$ be a compact suborbifold of
${\cal O}$. By definition, the preimage $\widetilde{\cal W}$ in
$\widetilde{\cal O}$ is a submanifold, and the space of orbifold
imbeddings $\Imb({\cal W},{\cal O})$ can be regarded as the quotient
of $\Imb_H(\widetilde{\cal W},\widetilde{\cal O})$ by the action of
$H$. For spaces of vectors, a subscript $H$ will indicate the
equivariant ones, thus for example ${\cal X}_H(\widetilde{\cal
W},T(\widetilde{\cal O}))$ means the $H$-equivariant sections of the
restriction of $T(\widetilde{\cal O})$ to $\widetilde{\cal W}$,
satisfying conditions~(1) and~(2) given in section~\ref{palais}.
Provided that $H$ acts as isometries on the $H$-invariant subset
$\widetilde{L}$ of $\widetilde{\cal O}$, the evaluation map $F$
carries ${\cal X}_H^L(T(\widetilde{\cal O}))$ into
$\Maps_H^{\widetilde{L}}(\widetilde{\cal O},\widetilde{\cal O})$.

The next two lemmas provide equivariant functions and metrics.

\begin{lemma}{equivariant function}{} Let $H$ be a group acting
smoothly and properly discontinuously on a manifold $M$, possibly with
boundary, such that $M/H$ is compact. Let $A$ be an $H$-invariant
closed subset of $M$, and $U$ an $H$-invariant neighborhood of
$A$. Then there exists an $H$-equivariant smooth function
$\gamma\colon M\to [0,1]$ which is identically equal to~$1$ on $A$ and
whose support is contained in~$U$.
\marginwrite{equivariant function}
\end{lemma}

\begin{proof}{} Fix a compact subset $C$ of $M$ which maps surjectively
onto $M/H$ under the quotient map.  Let $\phi\colon M\to[0,\infty)$ be
a smooth function such that $\phi(x)\geq 1$ for all $x\in C\cap A$ and
whose support is compact and contained in $U$. Define $\psi$ by
$\psi(x)\eq \sum_{h\in H}\phi(h(x))$. Now choose $\eta\colon\R\to[0,1]$
such that $\eta(r)\eq 0$ for $r\leq 0$ and $\eta(r)\eq 1$ for $r\geq
1$, and put $\gamma\eq \eta\circ\psi$.
\end{proof}

When ${\cal O}$ is compact, the following lemma provides a Riemannian
metric on $\widetilde{\cal O}$ for which the covering transformations
are isometries.

\begin{lemma}{covering isometries}{} Let $H$ be a group acting
smoothly and properly discontinuously on a manifold $M$, possibly with
boundary, such that $M/H$ is compact. Let $N$ be a properly imbedded
$H$-invariant submanifold, possibly empty. Then $M$ admits a complete
Riemannian metric, which is a product near $\partial M$ and such that
$N$ meets the collar $\partial M\times I$ in $I$-fibers, such that $H$
acts as isometries. Moreover, the action preserves the collar, and if
$(y,t)\in \partial M\times I$ and $h\in H$, then $h(y,t)\eq
(h\vert_{\partial M}(y),t)$.
\marginwrite{covering isometries}
\end{lemma}

\begin{proof}{} We first prove that equivariant Riemannian metrics
exist. Choose a compact subset $C$ of $M$ that maps surjectively onto
$M/H$ under the quotient map. Let $\phi\colon M\to [0,\infty)$ be a
compactly supported smooth function which is positive on $C$. Choose
a Riemannian metric $R$ on $M$ and denote by $R_x$ the inner product
which $R$ assigns to $T_x(M)$. Define a new metric $R'$ by
$$R'_x(v,w)=\sum_{h\in H}\phi(h(x))\,R_{h(x)}(h_*(v),h_*(w))\ .$$
\noindent Since $\phi$ is compactly supported, the sum is finite, and
since every orbit meets the support of $\phi$, $R'$ is positive
definite. To check equivariance, let $g\in H$. Then
$$\eqalign{R'_{g(x)}(g_*(v),g_*(w))&=
\sum_{h\in H}\phi(h(g(x)))\,
R_{h(g(x))}(h_*(g_*(v)),h_*(g_*(w)))\cr
&=\sum_{h\in H}\phi(hg(x))\,
R_{hg(x)}((hg)_*(v),(hg)_*(w))\cr
&=R'_x(v,w)\ .\cr}$$

\noindent We need to improve the metric near the boundary. First, note
that $C\cap\partial M$ maps surjectively onto the image of $\partial
M$. Choose an inward-pointing vector field $\tau'$ on a neighborhood
$U$ of $C\cap \partial M$, which is tangent to $N$. Choose a smooth
function $\phi\colon M\to [0,\infty)$ which is positive on $C\cap
\partial M$ and has support contained in $U$. The field $\phi\tau'$
defined on $U$ extends using the zero vector field on $M-U$ to a
vector field $\tau$ which is nonvanishing on $C\cap \partial M$. For
$x$ in the union of the $H$-translates of $U$, define
$\omega_x\eq\sum_{h\in H}\phi(h(x))\,h_*^{-1}(\tau_{h(x)})$. This is
defined, nonsingular, and equivariant on an equivariant neighborhood
of $\partial M$, and we use it to define a collar $\partial
M\times[0,2]$ equivariant in the sense that if $(y,t)\in\partial
M\times[0,2]$ then $h(y,t)\eq (h\vert_{\partial M}(y),t)$. Moreover,
$N$ meets this collar in $I$-fibers. On $\partial M\times[0,2]$,
choose an equivariant metric $R_1$ which is the product of an
equivariant metric on $\partial M$ and the standard metric on $[0,2]$,
and choose any equivariant metric $R_2$ defined on all of $M$. Using
lemma~\ref{equivariant function}, choose $H$-equivariant functions
$\phi_1$ and $\phi_2$ from $M$ to $[0,1]$ so that $\phi_1(x)\eq 1$ for
all $x\in \partial M\times [0,3/2]$ and the support of $\phi_1$ is
contained in $\partial M\times[0,2)$, and so that $\phi_2(x)\eq 1$ for
$x\in M-\partial M\times [0,3/2]$ and the support of $\phi_2$ is
contained in $M-\partial M\times[0,1]$.  Then, $\phi_1R_1+\phi_2R_2$
is $H$-equivariant and is a product near $\partial M$, and $N$ is
vertical in $\partial M\times I$.

Since $M/H$ is compact and $H$ acts as isometries, the metric must be
complete. For let $C$ be a compact subset of $M$ that maps
surjectively onto $M/H$. We may enlarge $C$ to a compact
codimension-zero submanifold $C'$ such that every point of $M$ has a
translate which lies in $C'$ at distance at least a fixed~$\epsilon$
from the frontier of $C'$. Then, any Cauchy sequence in $M$ can be
translated, except for finitely many terms, into a Cauchy sequence
in $C'$. Since $C'$ is compact, this converges, so the original sequence
also converged.
\end{proof}

We need the equivariant analogue of lemma~\ref{J}. Its proof uses
the following general fact.

\begin{proposition}{finitely generated}{} Suppose that $H$ acts
properly discontinuously on a locally compact connected space $X$, and
that $X/H$ is compact. Then $H$ is finitely generated.
\marginwrite{finitely generated}
\end{proposition}

\begin{proof}{} Using local compactness, there exists a compact set $C$
which maps surjectively to $X/H$. Let $H_0$ be the subgroup generated
by the finitely many elements $h$ such that $h(C)\cap C$ is nonempty.
The union of the $H_0$-translates of $C$ is an open and closed subset,
so must equal $X$. This implies that $H\eq H_0$.
\end{proof}

\begin{lemma}{orblemmaB}{}
Let $\widetilde{K}$ be an $H$-invariant subset of $\widetilde{\cal O}$
whose quotient in ${\cal O}$ is compact. Then there exists a
neighborhood $J$ of $1_{\widetilde{\cal O}}$ in
$\Maps_H^{\widetilde{K}}(\widetilde{\cal O}, \widetilde{\cal O})$ that
consists of diffeomorphisms.
\marginwrite{orblemmaB}
\end{lemma}

\begin{proof}{} Suppose first that ${\cal O}$ is compact. Then by
proposition~\ref{finitely generated}, $H$ is finitely generated. We
claim that if $f$ is a map that is close enough to $1_{\widetilde{\cal
O}}$, then $f$ commutes with the $H$-action. Choose an
$x\in\widetilde{\cal O}$ which is not fixed by any nontrivial element
of $H$.  Define $\Phi\colon \Maps_H^{\widetilde{K}}(\widetilde{\cal
O},\widetilde{\cal O})\to \hbox{End}(H)$ by sending $f$ to $\phi_f$
where $f(h(x))\eq\phi_f(h)f(x)$. This is independent of the choice of
$x$, hence is a homomorphism. If $f$ is close enough to
$1_{\widetilde{\cal O}}$ on $\set{x,h_1(x),\ldots,h_n(x)}$, where
$\set{h_1,\ldots,h_n}$ generates $H$, then $\phi_f\eq 1_H$. This prove
the claim.

Next we show that for $f$ close enough to $1_{\widetilde{\cal O}}$,
$f^{-1}(S)$ is compact whenever $S$ is compact.  From above, we may
assume that $f$ commutes with the $H$-action.  Let $C$ be a compact
set in $\widetilde{\cal O}$ which maps surjectively to ${\cal O}$.  If
$S$ is a set for which $f^{-1}(S)$ meets infinitely many translates of
$C$ then so does $S$, and $S$ could not be compact.

Consider $f$ close enough to $1_{\widetilde{\cal O}}$ to ensure the
previous conditions. By requiring $f$ sufficiently $C^\infty$-close to
$1_{\widetilde{\cal O}}$, $f_*$ is nonsingular at each point of $C$,
hence on all of $\widetilde{\cal O}$. If follows that $f$ is a local
diffeomorphism. Since also $f$ takes boundary to boundary and
preimages of compact sets are compact, $f$ is a covering map. Since
$\widetilde{\cal O}$ is simply-connected, $f$ is a diffeomorphism.

Now suppose that ${\cal O}$ is noncompact. Enlarge $\widetilde{K}/H$
to a compact codimen\-sion-zero suborbifold $L$. Let $\widetilde{L}'$
be a single component of $\widetilde{L}$ and $H'$ the stabilizer of
$\widetilde{L}'$ in $H$. Let $f\in\Maps_H^{\widetilde{K}}
(\widetilde{\cal O},\widetilde{\cal O})$. If $f$ is close enough to
$1_{\widetilde{\cal O}}$, then $f(\widetilde{L}')\eq\widetilde{L}'$
and by the previous argument we may assume that $f$ is a covering map
on $\widetilde{L}'$ (although since we don't know that
$\widetilde{L}'$ is simply connected, we cannot immediately conclude
that $f$ is a diffeomorphism). Since ${\cal O}$ is connected and
noncompact, $L$ has frontier in ${\cal O}$, hence $\widetilde{L}'$ has
frontier in $\widetilde{\cal O}$. Since $f$ is the identity on
$\widetilde{\cal O}-\widetilde{K}$, $f$ must be a diffeomorphism on
$\widetilde{L}'$, hence on all of $\widetilde{L}$, hence on all
of~$\widetilde{\cal O}$.
\end{proof}

We now prove the analogues of lemmas~\ref{extension}
and~\ref{logarithm} for vector fields on ${\cal O}$. Assume that
${\cal W}$ is a compact suborbifold of ${\cal O}$.

\begin{lemma}{orbextension}{} Let ${\cal W}$ be a compact suborbifold
of ${\cal O}$.
Let $L$ be a neighborhood of ${\cal W}$ in ${\cal O}$ and let $S$ be a
closed neighborhood in $\partial{\cal O}$ of $S\cap\partial{\cal W}$.
Denote the preimages in $\widetilde{\cal O}$ by $\widetilde{L}$ and
$\widetilde{S}$. Then there exists a continuous map $k\colon ({\cal
X}_H)_{<1/2}(\widetilde{\cal W},T(\widetilde{\cal O}))\to {\cal
X}_H^{\widetilde{L}}(T(\widetilde{\cal O}))$ such that $k(X)(x)\eq
X(x)$ for all $x$ in $\widetilde{\cal W}$ and $X\in({\cal
X}_H)_{<1/2}(\widetilde{\cal W},T(\widetilde{\cal O}))$. Moreover,
$k(Z)\eq Z$, and if $X(x)\eq Z(x)$ for all $x\in
\widetilde{S}\cap\partial\widetilde{\cal W}$, then $k(X)(x)\eq Z(x)$
for all $x\in\widetilde{S}$.
\marginwrite{orbextension}
\end{lemma}

\begin{proof}{} Assume first that ${\cal W}$ has positive codimension.
Replacing $L$ by a compact orbifold neighborhood of ${\cal W}$ and
using lemma~\ref{covering isometries}, we may assume that $H$ acts as
isometries on $\widetilde{O}$, that the metric is a product near
$\partial\widetilde{\cal O}$, and that $\widetilde{\cal W}$ meets the
collar $\partial\widetilde{\cal O}\times I$ in $I$-fibers. Let
$\nu(\widetilde{\cal W})$ be the normal bundle, regarded as a
subbundle of the restriction of $T(\widetilde{\cal O})$ to
${\widetilde{\cal W}}$. For $\epsilon>0$, let
$\nu_\epsilon(\widetilde{\cal W})$ be the subspace of all vectors of
length less than $\epsilon$. Since ${\cal W}$ is compact and $H$ acts
as isometries on $\widetilde{L}$, $\Exp$ imbeds
$\nu_\epsilon(\widetilde{\cal W})$ as a tubular neighborhood of
$\widetilde{\cal W}$ for sufficiently small $\epsilon$. By choosing
$\epsilon$ small enough, we may assume that
$\Exp(\nu_\epsilon(\widetilde{\cal W}))\subset \widetilde{L}$, that
the fibers at points in $\widetilde{S}$ map into $\widetilde{S}$, and
that the fibers at points in $\partial\widetilde{\cal
O}-\widetilde{S}$ map into $\partial\widetilde{\cal O}-\widetilde{S}$.

Now use lemma~\ref{equivariant function} to choose an $H$-equivariant
smooth function $\alpha\colon\widetilde{\cal O}\to[0,1]$ which is
identically equal to~1 on $\widetilde{\cal W}$ and has support in
$j(\nu_{\epsilon/2}(\widetilde{\cal W}))$. The extension $k(X)$ can
now be defined exactly as in lemma~\ref{extension}. Note that since
$H$ acts as isometries, the parallel translation function $P$ is
$H$-equivariant, and the $H$-equivariance of $k(X)$ follows easily.

The case when $W$ has codimension zero is similar, using
$\nu^+_\epsilon(\widetilde{\cal W})$ as in the proof of
lemma~\ref{extension}.
\end{proof}

\begin{lemma}{orblogarithm}{}
For all sufficiently small positive $\delta$, there exists a
continuous map $X\colon (\Imb_H)_{<\delta}(\widetilde{\cal
W},\widetilde{\cal O}) \to\allowbreak({\cal
X}_H)_{<1/2}(\widetilde{\cal W},T(\widetilde{\cal O}))$ such that
$\Exp(X(x))\eq j(x)$ for all $x\in \widetilde{\cal W}$ and $j\in
(\Imb_H)_{<\delta}(\widetilde{\cal W},\widetilde{\cal O})$, and
moreover if $j(x)\eq i_{\widetilde{\cal W}}(x)$ then $X(j)(x)\eq
Z(x)$.
\marginwrite{orblogarithm}
\end{lemma}

\begin{proof}{} Replacing ${\cal O}$ by a compact orbifold
neighborhood of ${\cal W}$ and using lemma~\ref{covering isometries},
we may assume that $H$ acts as isometries on $\widetilde{O}$, that the
metric is a product near $\partial\widetilde{O}$, and that
$\widetilde{\cal W}$ meets the collar $\partial\widetilde{\cal
O}\times I$ in $I$-fibers. Let $N_\epsilon(\widetilde{\cal W})$ be
defined exactly as in section~\ref{palais}. By compactness of ${\cal
W}$, there exists a positive $\epsilon$ such that for every
$x\in\widetilde{\cal W}$, $\Exp\colon N_\epsilon(\widetilde{\cal
W})\cap T_x(\widetilde{\cal O})\to \widetilde{\cal O}$ is a
diffeomorphism to a open neighborhood of $x$ in~$\widetilde{\cal O}$,
contained in $\widetilde{L}$.  The proof is then essentially the same
as the proof of lemma~\ref{logarithm}.
\end{proof}

The fundamental result is the analogue of theorem~B of \cite{P}.

\begin{theorem}{orbtheoremB}{} Let ${\cal W}$ be a compact suborbifold
of ${\cal O}$. Let $S$ be a closed neighborhood in $\partial{\cal O}$
of $S\cap\partial{\cal W}$, and let $L$ be a neighborhood of ${\cal
W}$ in ${\cal O}$. Then $\Imb^L({\cal W},{\cal O}\rel S)$ admits local
$\Diff^L({\cal O}\rel S)$ cross-sections.
\marginwrite{orbtheoremB}
\end{theorem}

\begin{proof}{} By proposition~\ref{inclusion}, it suffices to find a
local cross-section at the inclusion $i_{\cal W}$. Choose a compact
neighborhood $K$ of ${\cal W}$ with $K\subseteq L$. Using
lemmas~\ref{orblogarithm} and~\ref{orbextension}, there exist
continuous maps $X\colon
(\Imb_H)_{<\delta}(\widetilde{W},\widetilde{\cal O}) \to({\cal
X}_H)_{<1/2}(\widetilde{\cal W}, T(\widetilde{\cal O}))$ and $k\colon
({\cal X}_H)_{<1/2}(\widetilde{\cal W},T(\widetilde{\cal O}))
\allowbreak\to {\cal X}_H^{\widetilde{K}}(T(\widetilde{\cal O}))$. Let
$J$ be a neighborhood as in lemma~\ref{orblemmaB}. On a sufficiently
small neighborhood $\widetilde{U}$ of $i_{\widetilde{\cal W}}$, the
composition $\widetilde{\chi}\eq FkX$ is defined and has image in $J$.
Let $U$ be the imbeddings of ${\cal W}$ in ${\cal O}$ which admit a
lift to $\widetilde{U}$. By choosing $\widetilde{U}$ small enough, we
may ensure that the lift of an element of $U$ is unique. Define $\chi$
to be $\widetilde{\chi}$ applied to the lift of an element of $U$ to
$\widetilde{U}$, followed by the projection of
$\Diff_H^{\widetilde{K}}(\widetilde{\cal O})$ to $\Diff^{K}({\cal
O})$.

For elements in $U\cap \Imb^K({\cal W},{\cal O}\rel S)$, each lift to
$\widetilde{U}$ that is sufficiently close to $i_{\widetilde{\cal W}}$
must agree with $i_{\widetilde{\cal W}}$ on $\widetilde{S}$. So $U$
may be chosen small enough so that if $j\in U$ then
$\widetilde{j}\in\Imb(\widetilde{\cal W}),\widetilde{\cal
O}\rel\widetilde{S})$. Then, $X(\widetilde{j}(x))\eq Z(x)$ for all
$x\in \widetilde{S}$, so $k(X)(x)\eq Z(x)$ for all
$x\in\widetilde{S}$. It follows that $\chi(j)\in\Diff({\cal O}\rel
S)$.
\end{proof}

\begin{corollary}{orbcoro2}{} Let ${\cal W}$ be a compact suborbifold
of ${\cal O}$, which is either properly imbedded or codimension-zero.
Let $S$ be a closed neighborhood in $\partial{\cal O}$ of
$S\cap\partial{\cal W}$, and let $L$ be a neighborhood of ${\cal W}$
in ${\cal O}$. Then the restriction $\Diff^L({\cal O}\rel S)\to
\Imb^L({\cal W},{\cal O}\rel S)$ is locally trivial.
\marginwrite{orbcoro2}
\end{corollary}

\begin{corollary}{orbcoro3}{} Let ${\cal V}$ and ${\cal W}$ be
suborbifolds of ${\cal O}$, with ${\cal W}\subset
{\cal V}$. Assume that ${\cal W}$ compact, and is either properly
imbedded or codimension-zero. Let $S$ be a closed neighborhood in
$\partial{\cal O}$ of $S\cap\partial{\cal W}$, and let $L$ be a
neighborhood of ${\cal W}$ in ${\cal O}$. Then the restriction
$\Imb^L({\cal V},{\cal O}\rel S) \to \Imb^L({\cal W},{\cal
O}\rel S)$ is locally trivial.
\marginwrite{orbcoro3}
\end{corollary}

\section{Singular fiberings}
\marginwrite{sfiber}
\label{sfiber}

We will say that a continuous surjection $p\colon \Sigma\to {\cal O}$
of compact connected orbifolds is a {\it singular fibering} if there
exists a commutative diagram
$$\vbox{\halign{\hfil$#$\hfil&\hfil$#$\hfil&\hfil$#$\hfil\cr
\widetilde{\Sigma}&\mapright{\widetilde{p}}&\widetilde{{\cal O}}\cr
\mapdown{\sigma}&&\mapdown{\tau}\cr
\Sigma&\mapright{p}&{\cal O}\cr}}$$

\noindent in which
\begin{enumerate}
\item[{\rm(i)}] $\widetilde{\Sigma}$ and $\widetilde{\cal O}$ are
manifolds, and $\sigma$ and $\tau$ are regular orbifold coverings with
groups of covering transformations $G$ and $H$ respectively,
\item[{\rm(ii)}] $\widetilde{p}$ is surjective and locally trivial,
and
\item[{\rm(iii)}] the fibers of $p$ and $\widetilde{p}$ are
path-connected.
\end{enumerate}

The class of singular fiberings includes many Seifert fiberings, for
example all compact 3-dimensional Seifert manifolds $\Sigma$ except
the lens spaces with one or two exceptional orbits (see for
example~\cite{Scott}). For some of those lens spaces, ${\cal O}$ fails
to have an orbifold covering by a manifold. On the other hand, it is a
much larger class than Seifert fiberings, because no structure as a
homogeneous space is required on the fiber.

For mappings there is a complete analogy with the fibered case, where
now $\Diff_f(\Sigma)$ is by definition the quotient of the group of
fiber-preserving $G$-equivariant diffeomorphisms
$(\Diff_G)_f(\widetilde{\Sigma})$ by its normal subgroup $G$, and so
on.  A suborbifold $W$ of $\Sigma$ is called {\it vertical} if it is a
union of fibers. In this case the preimage $\widetilde{W}$ of $W$ in
$\widetilde{\Sigma}$ is a submanifold, and we can speak of
$\Imb_f(W,\Sigma)$ and $\Imb_v(W,\Sigma)$.

Following our usual notations, we put $\partial_v\Sigma\eq
p^{-1}(\partial{\cal O})$ and $\partial_vW\eq
W\cap\partial_v\Sigma$.

Since ${\cal O}$ is compact, lemma~\ref{covering isometries} shows
that a (complete) Riemannian metric on $\widetilde{\cal O}$ can be
chosen so that $H$ acts as isometries, and moreover so that the metric
on $\widetilde{\cal O}$ is a product near the boundary.  Next we will
sketch how to obtain a $G$-equivariant metric which is a product near
$\partial_h\widetilde{\Sigma}$ and near
$\partial_v\widetilde{\Sigma}$. If $\partial_h\widetilde{\Sigma}$ is
empty, we simply apply lemma~\ref{covering isometries}. Assume that
$\partial_h\widetilde{\Sigma}$ is nonemtpy. Construct a
$G$-equivariant collar of $\partial_h\widetilde{\Sigma}$, and use it
to obtain a $G$-equivariant metric such that the $I$-fibers of
$\partial_h\widetilde{\Sigma}\times I$ are vertical. If
$\partial_v\widetilde{\Sigma}$ is also nonempty, put
$Y\eq\partial_h\widetilde{\Sigma}\cap\partial_v\widetilde{\Sigma}$. We
will follow the construction in the last paragraph of
section~\ref{exponent}. Denote the collar of
$\partial_h\widetilde{\Sigma}$ by
$\partial_h\widetilde{\Sigma}\times[0,2]_1$. Assume that the metric on
$\partial_h\widetilde{\Sigma}$ was a product on a collar
$Y\times[0,2]_2$ of $Y$ in $\partial_h\widetilde{\Sigma}$. Next,
construct a $G$-equivariant collar $\partial_v\widetilde{\Sigma}\times
[0,2]_2$ of $\partial_v\widetilde{\Sigma}$ whose $[0,2]_2$-fiber at
each point of $Y\times [0,2]_1$ agrees with the $[0,2]_2$-fiber of the
collar of $Y$ in $\partial_h\widetilde{\Sigma}\times\set{t}$. Then,
the product metric on $\partial_v\widetilde{\Sigma}\times[0,2]_2$
agrees with the product metric of
$\partial_h\widetilde{\Sigma}\times[0,2]_1$ where they overlap, and
the $G$-equivariant patching can be done to obtain a metric which is a
product near $\partial_v\widetilde{\Sigma}$ without losing the
property that it is a product near~$\partial_h\widetilde{\Sigma}$.  We
will always assume that the metrics have been selected with these
properties. In particular, $G$ preserves the vertical and horizontal
parts of vectors.

Some basic observations about singular fiberings will be needed.

\begin{lemma}{lift}{} The action of $G$ preserves the fibers of
$\widetilde{p}$. Moreover:
\begin{enumerate}
\item[\rm(i)] If $g\in G$, then there exists an element $h\in H$ such
that $\widetilde{p}g\eq h\widetilde{p}$.
\item[\rm(ii)] If $h\in H$, then there exists an element
$g$ of $G$ such that $\widetilde{p}g\eq h\widetilde{p}$.
\item[\rm(iii)] If $x\in \Sigma$, then $\tau^{-1}p(x)\eq
\widetilde{p}\sigma^{-1}(x)$.
\end{enumerate}
\marginwrite{lift}
\end{lemma}

\begin{proof}{}
Suppose that $\widetilde{p}(x)\eq\widetilde{p}(y)$. For $g\in G$, we
have $\tau\widetilde{p}(g(x))\eq p\sigma(g(x))\eq p\sigma(x)\eq
\tau\widetilde{p}(x)\eq \tau\widetilde{p}(y)\eq
\tau\widetilde{p}(g(y))$. Since the fibers of $\widetilde{p}$ are
path-connected, and the fibers of $\tau$ are discrete, this implies
that $g(x)$ and $g(y)$ lie in the same fiber of $\widetilde{p}$. For
(i), let $g\in G$. Since $g$ preserves the fibers of $\widetilde{p}$,
it induces a map $h$ on $\widetilde{{\cal O}}$. Given
$x\in\widetilde{{\cal O}}$, choose $y\in\widetilde{\Sigma}$ with
$\widetilde{p}(y)\eq x$. Then $\tau h(x)\eq \tau
\widetilde{p}(g(y))\eq p\sigma(g(y))\eq p\sigma(y)\eq \tau
\widetilde{p}(y)\eq \tau(x)$ so $h\in H$.

To prove (ii), suppose $h$ is any element of $H$. Let
$\hbox{sing}({\cal O})$ denote the singular set of ${\cal O}$. Choose
$a\in\widetilde{{\cal O}}-\tau^{-1}(\hbox{sing}({\cal O}))$, choose
$s\in \widetilde{\Sigma}$ with $\widetilde{p}(s)\eq a$, and choose
$s''\in \widetilde{\Sigma}$ with $\widetilde{p}(s'')\eq h(a)$. Since
$p\sigma(s)\eq \tau\widetilde{p}(s)\eq \tau\widetilde{p}(s'')\eq
p\sigma(s'')$, $\sigma(s)$ and $\sigma(s'')$ must lie in the same
fiber of $p$. Since the fiber is path-connected, there exists a path
$\beta$ in that fiber from $\sigma(s'')$ to $\sigma(s)$. Let
$\widetilde{\beta}$ be its lift in $\widetilde{\Sigma}$ starting at
$s''$ and let $s'$ be the endpoint of this lift, so that
$\sigma(s')\eq \sigma(s)$. Note that $\widetilde{p}(s')\eq
\widetilde{p}(s'')\eq h(a)$ since $\widetilde{\beta}$ lies in a fiber
of $\widetilde{p}$. Since $\sigma(s)\eq \sigma(s')$, there exists a
covering transformation $g\in G$ with $g(s)\eq s'$. To show that
$\widetilde{p}g\eq h\widetilde{p}$, it is enough to verify that they
agree on the dense set $\widetilde{p}^{-1}(\widetilde{{\cal
O}}-\tau^{-1}(\hbox{sing}({\cal O})))$. Let $t\in
\widetilde{p}^{-1}(\widetilde{{\cal O}}-\tau^{-1}(\hbox{sing}({\cal
O})))$ and choose a path $\gamma$ in
$\widetilde{p}^{-1}(\widetilde{{\cal O}}-\tau^{-1}(\hbox{sing}({\cal
O})))$ from $s$ to $t$. Since $g\in G$, we have $p\sigma\gamma\eq
p\sigma g\gamma$. Therefore $\tau\widetilde{p}\gamma\eq
\tau\widetilde{p} g\gamma$, and so $\widetilde{p}g\gamma$ is the
unique lift of $p\sigma\gamma$ starting at $\widetilde{p}g(s)\eq
h(a)$. But this lift equals $h\widetilde{p}\gamma$, so
$h\widetilde{p}(t)\eq \widetilde{p}g(t)$.

For (iii), fix $z_0\in \sigma^{-1}(x)$ and let $y_0\eq
\widetilde{p}(z_0)$. Suppose $y\in \widetilde{p}\sigma^{-1}(x)$.
Choose $z\in \sigma^{-1}(x)$ with $\widetilde{p}(z)\eq y$. Since
$\sigma$ is a regular covering, there exists $g\in G$ such that $g(z)\eq
z_0$. By~(i), $g$ induces $h$ on $\widetilde{\cal O}$, and $h(y)\eq
h\widetilde{p}(z)\eq \widetilde{p}g(z)\eq \widetilde{p}(z_0)\eq y_0$.
Therefore $\tau(y)\eq \tau(h(y))\eq \tau(y_0)\eq
\tau\widetilde{p}(z_0)\eq p\sigma(z_0)\eq p(x)$ so $y\in
\tau^{-1}(p(x))$. For the opposite inclusion, suppose that $y\in
\tau^{-1}p(x)$, so $\tau(y)\eq p(x)\eq \tau(y_0)$. Since $\sigma$ is
regular, there exists $h\in H$ such that $h(y_0)\eq y$. Let $g$ be as
in (ii). Then $y\eq h(y_0)\eq h\widetilde{p}(z_0)\eq
\widetilde{p}g(z_0)$, and $\sigma(g(z_0))\eq \sigma(z_0)\eq x$ so
$y\in\widetilde{p}(\sigma^{-1}(x))$.
\end{proof}

One consequence of lemma~\ref{lift} is that (smooth nonsingular) paths
in $\widetilde{\cal O}$ have horizontal lifts in $\widetilde{\Sigma}$.
To see this, we first claim that the horizontal lifts of any vector
$\omega$ in $T(\widetilde{\cal O})$ have bounded lengths. Fix a compact
subset $C$ of $\widetilde{\Sigma}$ such that $\sigma(C)\eq\Sigma$.
Let $H_*\omega$ be the set of $H$-translates of $\omega$. Since $C$ is
compact and $H_*\omega$ is closed, the lengths of the horizontal lifts
of vectors in $H_*\omega$ to vectors in $T(\widetilde{\Sigma})\vert_C$
are bounded by some~$L$. If $\widetilde{\omega}$ is any lift of
$\omega$, there exists $g\in G$ such that $g_*\widetilde{\omega}\in
T(\widetilde{\Sigma})\vert_C$. By lemma~\ref{lift}(i), there exists
$h\in H$ such that $\widetilde{p}_*(g_*\widetilde{\omega})\eq
h_*\omega$. Since $g_*\widetilde{\omega}$ is a horizontal lift of
$h_*(\omega)$ and $g_*$ is an isometry, $\norm{\widetilde{\omega}}\eq
\norm{g_*\widetilde{\omega}}\leq L$, proving the claim. Since the
metric on $\widetilde{\Sigma}$ is complete, the claim shows that a
path in $\widetilde{\cal O}$ could only fail to lift if a partial lift
started in $\widetilde{\Sigma}-\partial_h\widetilde{\Sigma}$ and then
reached a point of $\partial_h\widetilde{\Sigma}$, impossible since
the metric is a product near~$\partial_h\widetilde{\Sigma}$.

Since horizontal lifts exist, the aligned exponential $\Exp_a$ of
$\widetilde{\Sigma}$ is defined. Moreover, it is $G$-equivariant:
since $G$ consists of fiber-preserving isometries, $\Exp_v$ is
$G$-equivariant, and since $G$ preserves horizontal parts of vectors,
it preserves horizontal lifts.

The notations ${\cal A}(\widetilde{W},T(\widetilde{\Sigma}))$ and
${\cal A}(T(\widetilde{\Sigma}))$ and the map $F_a\colon {\cal
A}(T(\widetilde{\Sigma}))\to
\Maps(\widetilde{\Sigma},\widetilde{\Sigma})$ are analogous to those
in section~\ref{project}.

\begin{theorem}{sftheorem1}{} Let $S$ be a closed subset of
$\partial{\cal O}$, and let $T\eq p^{-1}(S)$. Then
$\Diff({\cal O}\rel S)$ admits local
$\Diff_f(\Sigma\rel T)$ cross-sections.
\marginwrite{sftheorem1}
\end{theorem}

\begin{proof}{} Lemma~\ref{orblogarithm}, with
${\cal W}\eq{\cal O}$, provides $X\colon
\widetilde{U}_\delta
\to ({\cal X}_H)_{<1/2}(T(\widetilde{\cal O}))$, where
$\widetilde{U}_\delta\eq\set{f\in \Diff_H(\widetilde{\cal O})\vbar
d(f(x),x)<\delta\hbox{\ for all $x\in \widetilde{\cal O}$}}$.
Let $h\in \Diff({\cal O})$ and let $\widetilde{h}\in
\Diff_H(\widetilde{\cal O})$ be a lift of $h$.
For every $\widetilde{g}\in \widetilde{U}_\delta \widetilde{h}$,
$\Exp(X(\widetilde{g}\widetilde{h}^{-1}(x)))\eq
\widetilde{g}\widetilde{h}^{-1}(x)$. Define
$\widetilde{\chi}\colon\widetilde{U}_\delta\to {\cal
A}_G(T(\widetilde{\Sigma}))$ by
$$\widetilde{\chi}(\widetilde{g})(x)=
\big(\widetilde{p}_*\vert_{H_x}\big)^{-1}(X(\widetilde{g}
\widetilde{h}^{-1})(\widetilde{p}(x)))\ .$$

\noindent The boundary tangency conditions are clearly satisfied, and
$\Exp_a(\widetilde{\chi}(\widetilde{g})(x))$ exists since it is the
horizontal lift of a geodesic from $\widetilde{p}(x)$ to
$\widetilde{g}\widetilde{h}^{-1}(\widetilde{p}(x))$.  To see that
$\widetilde{\chi}(\widetilde{g})$ is $G$-equivariant, suppose
$\gamma\in G$. By lemma~\ref{lift}(i), $\gamma$ induces $\lambda\in
H$. We have
$$\eqalign{\widetilde{\chi}(\widetilde{g})(\gamma(x))
&=\big(\widetilde{p}_*\vert_{H_x}\big)^{-1}(X(\widetilde{g}
\widetilde{h}^{-1})(\widetilde{p}(\gamma(x))))\cr
&=\big(\widetilde{p}_*\vert_{H_x}\big)^{-1}(X(\widetilde{g}
\widetilde{h}^{-1})(\lambda\widetilde{p}(x)))\cr
&=\big(\widetilde{p}_*\vert_{H_x}\big)^{-1}(\lambda_*X(\widetilde{g}
\widetilde{h}^{-1})(\widetilde{p}(x)))\cr
&=\gamma_*\big(\widetilde{p}_*\vert_{H_x}\big)^{-1}(X(\widetilde{g}
\widetilde{h}^{-1})(\widetilde{p}(x)))\cr
&=\gamma_*\widetilde{\chi}(\widetilde{g})(x)\ ,\cr}$$

\noindent the penultimate equality using the fact that $G$ preserves the
horizontal subspaces.

Let $\widetilde{U}\eq \widetilde{\chi}^{-1}(J)$, where $J$ is a
neighborhood of $1_{\widetilde{\Sigma}}$ as in lemma~\ref{orblemmaB}.
Let $U$ be a neighborhood of $h$ consisting of elements having a lift
in $\widetilde{U}$. Since $G$ is a discrete subgroup of
$\Diff_G(\widetilde{\Sigma})$, we may choose $\delta$ small enough to
ensure that these lifts are unique. Now we can define $\chi\colon U\to
\Diff_f(\Sigma)$ by putting $\chi(g)$ equal to the diffeomorphism
induced on $\Sigma$ by $F_a\widetilde{\chi}(\widetilde{g})$.
\end{proof}

From proposition~\ref{theoremA}, we have immediately

\begin{theorem}{sfproject diffs}{} Let $S$ be a closed subset of
$\partial{\cal O}$, and let $T\eq p^{-1}(S)$. Then
$\Diff_f(\Sigma\rel T)\to \Diff({\cal O}\rel S)$ is
locally trivial.
\marginwrite{sfproject diffs}
\end{theorem}

We now extend lemmas~\ref{lemmaC} and~\ref{lemmaD} to the singular
fibered case.

\begin{lemma}{sflemmaC}{}
Let $W$ be a vertical suborbifold of $\Sigma$. Let $T$ be a closed
fibered neighborhood in $\partial_v\Sigma$ of $T\cap \partial_vW$.
Then for all sufficiently small $\delta$, there exists a continuous
map $k\colon({\cal
A}_G)_{<\delta}(\widetilde{W},T(\widetilde{\Sigma}))\to {\cal
A}_G(T(\widetilde{\Sigma}))$ such that $k(X)(x)\eq X(x)$ for all $x\in
\widetilde{W}$ and $X\in({\cal
A}_G)_{<\delta}(\widetilde{W},T(\widetilde{\Sigma}))$. If $X(x)\eq
Z(x)$ for all $x\in \widetilde{T}\cap\partial_v\widetilde{W}$, then
$k(X)(x)\eq Z(x)$ for all $x\in \widetilde{T}$. Furthermore, $k(({\cal
V}_G)_{<\delta}(\widetilde{W}, T(\widetilde{\Sigma})))\subset{\cal
V}_G(T(\widetilde{\Sigma}))$.
\marginwrite{sflemmaC}
\end{lemma}

\begin{proof}{} As with lemma~\ref{extension}, the positive
codimension and codimen\-sion-zero cases are similar, so we only discuss
the former. Let ${\cal W}$ be the image of $W$ in ${\cal O}$, and
denote $\tau^{-1}{\cal W}$ by $\widetilde{\cal W}$. By
lemma~\ref{lift}(iii), $\widetilde{\cal W}\eq
\widetilde{p}(\widetilde{W})$, and by lemma~\ref{lift}(ii), it is
$H$-invariant. Since it is a submanifold of $\widetilde{\cal O}$, it
follows that ${\cal W}$ is a suborbifold of ${\cal O}$. A section
$X\in{\cal A}_G(\widetilde{W},T(\widetilde{\Sigma}))$ induces a
well-defined section $\widetilde{p}_*X\in {\cal X}(\widetilde{\cal
W},T(\widetilde{\cal O}))$. By lemma~\ref{lift}(ii),
$\widetilde{p}_*X$ is $H$-equivariant.

We claim that there exists a positive $\delta$ so that if $X\in ({\cal
A}_G)_{<\delta}(\widetilde{W},T(\widetilde{\Sigma}))$ then
$\widetilde{p}_*X\in ({\cal X}_H)_{<1/2}(\widetilde{\cal
W},T(\widetilde{\cal O}))$. For if not, there would be a sequence $x_i$
in $\widetilde{W}$ such that $\norm{X(x_i)}\to 0$ but
$\norm{\widetilde{p}_*X(\widetilde{p}(x_i))}\geq 1/2$.  Since $W$ is
compact, there exists a compact subset $C\subset
\widetilde{W}$ such that $\sigma(C)\eq W$. There exist
elements $g_i\in G$ so that $g_i(x_i)\in C$, and if $h_i\in H$ are
obtained using lemma~\ref{lift}(i) then
$\norm{X(g_i(x_i))}\eq\norm{X(x_i)}$ while
$\norm{\widetilde{p}_*X(\widetilde{p}g_i(x_i))}\eq
\norm{\widetilde{p}_*X(h_i\widetilde{p}(x_i))}\eq
\norm{\widetilde{p}_*X(\widetilde{p}(x_i))}\allowbreak\geq 1/2$. So
we may assume that the $x_i$ lie in $C$, hence that they converge to
$x\in C$. Then, $\norm{X(x)}\eq 0$ but
$\norm{\widetilde{p}_*X(\widetilde{p}(x))}\geq 1/2$, a contradiction.

We now follow the proof of lemma~\ref{lemmaC}.  Let
$k_{\widetilde{\cal O}}\colon ({\cal X}_H)_{<1/2}(\widetilde{\cal
W},T(\widetilde{\cal O}))\to {\cal X}_H(T(\widetilde{\cal O}))$ be
obtained using lemma~\ref{orbextension}.  Let
$\nu_\epsilon(\widetilde{W})$ be the $\epsilon$-normal bundle of
$\widetilde{W}$. Since $W$ is compact, for sufficiently
small~$\epsilon$, $j_a\colon \nu_\epsilon(\widetilde{W})\to
\widetilde{\Sigma}$ defined by
$j_a(\omega)\eq\Exp_a(\omega)$ and carries
$\nu_\epsilon(\widetilde{W})$ diffeomorphically to a neighborhood of
$\widetilde{W}$ in $\widetilde{\Sigma}$.  Since $W$ is compact, we may
choose $\epsilon$ small enough so that $j_a(\omega)\in
\partial_h\widetilde{\Sigma}\times I$ only
when $\pi(\omega)\in\partial_h\widetilde{\Sigma}\times I$.

Since $G$ acts as isometries and preserves horizontal lifts, the
aligned parallel translation $P_a$ is $G$-equivariant. Using
lemma~\ref{equivariant function} there exists a smooth $G$-equivariant
function $\alpha\colon \widetilde{\Sigma}\to [0,1]$ which is
identically~1 on $\widetilde{W}$ and identically~0 on
$\widetilde{\Sigma}-j(\nu_{\epsilon/2}(\widetilde{W}))$.  Define
$k_{\widetilde{\Sigma}}\colon ({\cal
A}_G)_{<\delta}(\widetilde{W},T(\widetilde{\Sigma}))\to {\cal
V}_G(T(\widetilde{\Sigma}))$ by
$$k_{\widetilde{\Sigma}}(X)(x)\eq
\cases{
\alpha(x)P_a(X(\pi(j_a^{-1}(x))),j_a^{-1}(x))_v&for
$x\in j_a(\nu_{\epsilon}(\widetilde{W}))$\cr
Z(x)&for $x\in \widetilde{\Sigma}-j_a(\nu_{\epsilon/2}(\widetilde{W}))$\cr}$$

\noindent and $k$ by
$$k(X)(x)=(\widetilde{p}_*\vert_{H_x})^{-1}
(k_{\widetilde{\cal O}}(p_*X)(p(x)))
\;+\;k_{\widetilde{\Sigma}}(X_v)(x)\ .$$
\end{proof}

\begin{lemma}{sflemmaD}{} Let $W$ be a vertical suborbifold of
$\Sigma$. For small $\delta>0$, there exists a continuous map
$$X\colon ((\Imb_G)_f)_{<\delta}(\widetilde{W},
\widetilde{\Sigma})
\to{\cal A}_G(\widetilde{W},
T(\widetilde{\Sigma}))$$
such that $\Exp_a(X(j)(x))\eq j(x)$ for all
$x\in \widetilde{W}$ and $j\in ((\Imb_G)_f)_{<\delta}(\widetilde{W},
\widetilde{\Sigma})$.
Moreover, $X(((\Imb_G)_v)_{<\delta}(\widetilde{W},
\widetilde{\Sigma}))\subseteq{\cal
V}_G(\widetilde{W},T(\widetilde{\Sigma}))$, and if $j(x)\eq
i_{\widetilde{W}}(x)$ then $X(j)(x)\eq Z(x)$.
\marginwrite{sflemmaD}
\end{lemma}

\begin{proof}{} Let $N_\epsilon(\widetilde{W})$ be as defined before
the proof of lemma~\ref{lemmaD}. Since $W$ is compact, we can choose
$\epsilon$ small enough to ensure the local diffeomorphism
condition. Choose $\delta$ small enough so that
$j(x)\in\Exp_a(N_\epsilon(\widetilde{W})\cap T_x(\widetilde{\Sigma}))$
for every $x\in \widetilde{W}$ and $j\in
((\Imb_G)_f)_{<\delta}(\widetilde{W}, \widetilde{\Sigma})$. Define
$X(j)(x)$ to be the unique vector in $N_\delta(\widetilde{W})\cap
T_x(\widetilde{\Sigma})$ such that $\Exp_a(X(j)(x))$ equals~$j(x)$.
\end{proof}

\begin{theorem}{sftheorem2}{}
Let $W$ be a vertical suborbifold of $\Sigma$. Let $T$ be a closed
fibered neighborhood in $\partial_v\Sigma$ of $T\cap \partial_vW$.
Then
\begin{enumerate}
\item[{\rm (i)}] $\Imb_f(W,\Sigma\rel T)$ admits local
$\Diff_f(\Sigma\rel T)$ cross-sections, and
\item[{\rm (ii)}] $\Imb_v(W,\Sigma\rel T)$ admits local
$\Diff_v(\Sigma\rel T)$ cross-sections.
\end{enumerate}
\marginwrite{sftheorem2}
\end{theorem}

\begin{proof}{} By proposition~\ref{inclusion}, it suffices to
construct local cross-sections at the inclusion $i_W$. Obtain
$k\colon({\cal A}_G)_{<\delta}(\widetilde{W},T(\widetilde{\Sigma}))\to
({\cal A}_G)_{<1/2}(T(\widetilde{\Sigma}))$, and $X\colon
((\Imb_G)_f)_{<\delta_1}(\widetilde{W},\widetilde{\Sigma}) \to {\cal
A}_G(\widetilde{W},T(\widetilde{\Sigma}))$ using lemmas~\ref{sflemmaC}
and~\ref{sflemmaD}. Fix a neighborhood $\widetilde{U}$ of
$i_{\widetilde{W}}$ small enough so that $X(\widetilde{U})\subseteq
({\cal A}_G)_{<\delta_1}(\widetilde{W}, T(\widetilde{\Sigma}))$. Let
$U$ be a neighborhood of $i$ small enough so that each element $j$ of
$U$ has a unique lift $\widetilde{j}$ into $\widetilde{U}$, and so
that if $j$ agrees with $i_W$ on $\partial_vW$ then $\widetilde{j}$
agrees with $i_{\widetilde{W}}$ on $\partial_v\widetilde{W}$. For
$j\in U$, define $\chi(j)$ to be the element of $\Diff_f(\Sigma\rel
T)$ induced by $F_akX(\widetilde{j})$.
\end{proof}

As in section~\ref{restrict}, we have the following immediate
corollaries.

\begin{corollary}{sfcorollary2}{} Let $W$ be a vertical suborbifold of
$\Sigma$. Let $T$ be a fibered neighborhood in $\partial_v\Sigma$ of
$T\cap\partial_vW$. Then the following restrictions are locally
trivial:
\begin{enumerate}
\item[{\rm(i)}] $\Diff_f(\Sigma\rel T)\to \Imb_f(W,\Sigma\rel T)$, and
\item[{\rm(ii)}] $\Diff_v(\Sigma\rel T)\to
\Imb_v(W,\Sigma\rel T)$.
\end{enumerate}
\marginwrite{sfcorollary2}
\end{corollary}

\begin{corollary}{sfcorollary3}{} Let $V$ and $W$ be vertical
suborbifolds of $\Sigma$, with $W\subseteq V$. Let $T$ be a closed
fibered neighborhood in $\partial_v\Sigma$ of $T\cap\partial_vW$. Then
the following restrictions are locally trivial:
\begin{enumerate}
\item[{\rm(i)}] $\Imb_f(V,\Sigma\rel T)\to
\Imb_f(V,\Sigma\rel T)$, and
\item[{\rm(ii)}] $\Imb_v(V,\Sigma\rel T)\to
\Imb_v(W,\Sigma\rel W)$.
\end{enumerate}
\marginwrite{sfcorollary3}
\end{corollary}

\begin{theorem}{sfsquare}{} Let $W$ be a vertical suborbifold of
$\Sigma$. Let $T$ be a closed fibered neighborhood in
$\partial_v\Sigma$ of $T\cap\partial_vW$, and let $S\eq p(T)$. Then
all four maps in the following square are locally trivial:
$$\vbox{\halign{\hfil#\hfil\quad&#&\quad\hfil#\hfil\cr
$\Diff_f(\Sigma\rel T)$&$\longrightarrow$&%
$\Imb_f(W,\Sigma\rel T)$\cr
\noalign{\smallskip}
$\mapdown{}$&&$\mapdown{}$\cr
\noalign{\smallskip}
$\Diff({\cal O}\rel S)$&$\longrightarrow$&%
$\Imb(p(W),{\cal O}\rel S)$\rlap{\ .}\cr}}$$
\marginwrite{sfsquare}
\end{theorem}

\section{Restricting to the boundary or the basepoint}
\label{basept}
\marginwrite{basept}

Our restriction theorems deal with the case when the suborbifold is
properly imbedded. By a simple doubling trick, we can also extend to
restriction to suborbifolds of the boundary.

\begin{proposition}{restrict to boundary}{} Let $\Sigma\to{\cal O}$ be
a singular fibering. Let $S$ be a suborbifold of $\partial{\cal O}$,
and let $T\eq p^{-1}(S)$. Then
\begin{enumerate}
\item[{\rm(a)}] $\Imb(S,\partial{\cal O})$ admits local
$\Diff({\cal O})$ cross-sections.
\item[{\rm(b)}] $\Imb_f(T,\partial_v{\Sigma})$ admits local
$\Diff_f(\Sigma)$ cross-sections.
\end{enumerate}
\marginwrite{restrict to boundary}
\end{proposition}

\begin{proof}{} For (a), we first show that $\Diff(\partial{\cal O})$
admits local $\Diff({\cal O})$ cross-sections. Let $\Delta$ be the
double of ${\cal O}$ along $\partial{\cal O}$, and regard ${\cal O}$ as
a suborbifold of $\Delta$ by identifying it with one of the two copies
of ${\cal O}$ in $\Delta$. By theorem~\ref{orbtheoremB},
$\Imb(\partial{\cal O}, \Delta)$ admits local $\Diff(\Delta)$
cross-sections. We may regard $\Diff(\partial{\cal O})$ as a subspace
of $\Imb(\partial{\cal O},\Delta)$. Suppose $\chi\colon
U\to\Diff(\Delta)$ is a local cross-section at a point in
$\Imb(\partial{\cal O},\Delta)$ that lies in $\Diff(\partial{\cal
O})$. By composing with the diffeomorphism of $\Delta$ that
interchanges the two copies of ${\cal O}$, and reducing the size of
$U$ if necessary, we may assume that $\chi$ carries the elements of
$U$ that preserve $\partial{\cal O}$ to diffeomorphisms that preserve
${\cal O}$. Then a local $\Diff({\cal O})$ cross-section on
$U\cap\Diff(\partial{\cal O})$ is defined by sending $g$ to
$\chi(g)\vert_{\cal O}$.

By proposition~\ref{inclusion}, for (a) it suffices to produce local
cross-sections at the inclusion $i_S$. By theorem~\ref{orbtheoremB},
there is a local $\Diff(\partial{\cal O})$ cross-section $\chi_1$ for
$\Imb(S,\partial {\cal O})$ at $i_S$. Let $\chi_2$ be a local
$\Diff({\cal O})$ cross-section for $\Diff(\partial {\cal O})$ at
$\chi_1(i_S)$. On a neighborhood $U$ of $i_S$ in
$\Imb(S,\partial{\cal O})$ small enough so that $\chi_2\chi_1$ is
defined, the composition is the desired $\Diff({\cal O})$
cross-section. For if $j\in U$, then
$\chi_2(\chi_1(j))(i_S)(x)=\chi_1(j)(i_S)(x)=\chi_1(j)(x)=j(x)$.

The proof of (b) is similar. Double $\Sigma$ along $\partial_v\Sigma$
and apply theorem~\ref{sftheorem2}, to produce local $\Diff_f(\Sigma)$
cross-sections for $\Diff_f(\partial_v\Sigma)$. Apply it again to
produce local $\Diff_f(\partial_v\Sigma)$ cross-sections for
$\Imb_f(T,\partial_v\Sigma)$. Their composition, where defined, is a
local $\Diff_f(\Sigma)$ cross-section for
$\Imb_f(T,\partial_v\Sigma)$.
\end{proof}

An immediate consequence is

\begin{corollary}{special1}{} Let $\Sigma\to{\cal O}$ be a singular
fibering. Let $S$ be a suborbifold of $\partial{\cal O}$, and let
$T=p^{-1}(S)$. Then $\Diff({\cal O})\to \Imb(S,\partial{\cal O})$ and
$\Diff_f(\Sigma)\to \Imb_f(T,\partial_v\Sigma)$ are locally trivial.
In particular, $\Diff({\cal O})\to\Diff(\partial{\cal O})$ and
$\Diff_f(\Sigma)\to\Diff_f(\partial_v\Sigma)$ are locally trivial.
\marginwrite{special1}
\end{corollary}

Here are two other consequences which are applied in \cite{M-R}.

\begin{corollary}{special2}{} Let ${\cal W}$ be a suborbifold of ${\cal
O}$. Then $\Imb({\cal W},{\cal O})\to\Imb({\cal W}\cap \partial{\cal
O},\partial{\cal O})$ is locally trivial.
\marginwrite{special2}
\end{corollary}

\begin{proof}{} By theorem~\ref{orbtheoremB}, $\Imb({\cal W}\cap
\partial{\cal O},\partial{\cal O})$ admits local $\Diff(\partial{\cal
O})$ cross-sections, and by proposition~\ref{restrict to boundary},
$\Diff(\partial{\cal O})$ admits local $\Diff({\cal O})$
cross-sections. Composing them gives local $\Diff({\cal O})$
cross-sections for $\Imb({\cal W}\cap\partial{\cal O},\partial{\cal
O})$.
\end{proof}

\begin{corollary}{special3}{} Let $W$ be a vertical suborbifold of
$\Sigma$. Then $\Imb_f(W,\Sigma)\to \Imb_f(W\cap
\partial_v\Sigma,\partial_v\Sigma)$ is locally trivial.
\marginwrite{special3}
\end{corollary}

\begin{proof}{} Theorem~\ref{sftheorem2}, applied to
$\partial_v\Sigma$, and proposition~\ref{restrict to boundary} show
that $\Imb_f(W\cap\partial_v\Sigma,\partial_v\Sigma)$ admits local
$\Diff_f(\Sigma)$ cross-sections.
\end{proof}

Many applications of the fibration $\Diff(M)\to\Imb(V,M)$ concern the
case when the submanifold is a single point. Since in the fibered case
a single point is not usually a vertical submanifold, this case is not
directly covered by our previous theorems. The next proposition allows
nonvertical suborbifolds that are contained in a single fiber, so
applies when the submanifold is a single point. To set notation, let
$p\colon \Sigma\to {\cal O}$ be a singular fibering. Let $P$ be a
suborbifold of $\Sigma$ which is contained in a single fiber~$F$. Let
$T$ be a fibered closed subset of $\partial_v\Sigma$. By
$\Imb_t(P,\Sigma\rel T)$ we denote the orbifold imbeddings whose image
is contained in a single fiber of~$\Sigma$, which restrict to the
identity on $P\cap T$, and which map $P\cap(\partial \Sigma-T)$ into
$\partial\Sigma-T$.

\begin{proposition}{restrict to basepoint}{} Let $T$ be a fibered
closed subset of $\partial_v\Sigma$, which is a neighborhood in
$\partial_v\Sigma$ of $P\cap T$. Then $\Imb_t(P,\Sigma\rel T)$ admits
local $\Diff_f(\Sigma\rel T)$ cross sections.
\marginwrite{restrict to basepoint}
\end{proposition}

\begin{proof}{}
Notice that $p(P)$ is a point and is a properly imbedded suborbifold
of ${\cal O}$, with orbifold structure determined by the local group
at $p(P)$. Each imbedding $i\in\Imb_t(P,\Sigma)$ induces an orbifold
imbedding $pi\colon p(P)\to{\cal O}$. Let $S\eq p(T)$.

By proposition~\ref{inclusion}, it suffices to produce a local
cross-section at the inclusion $i_P$. By theorem~\ref{orbtheoremB},
$\Imb(p(P),{\cal O}\rel S)$ has local $\Diff({\cal O}\rel S)$
cross-sections, and by proposition~\ref{sftheorem1}, $\Diff({\cal
O}\rel S)$ has local $\Diff_f(\Sigma\rel T)$ cross-sections. A
suitable composition of these gives a local $\Diff_f(\Sigma\rel T)$
cross-section $\chi_1$ for $\Imb(p(P),{\cal O}\rel S)$ at $pi_P$. As
remarked in section~\ref{palais}, we may assume that $\chi_1(pi_P)$ is
the identity diffeomorphism of $\Sigma$. By corollary~\ref{orbcoro2},
there exists a local $\Diff(F\rel T\cap F)$ cross-section $\chi_2$ for
$\Imb(P,F\rel T\cap F)$ at $i_P$, and we may assume that $\chi_2(i_P)$
is the identity diffeomorphism of $F$. Let $\chi_3$ be a local
$\Diff_f(\Sigma\rel T)$ cross-section for $\Imb_f(F,\Sigma\rel T)$ at
$i_F$ given by corollary~\ref{sfcorollary2}. Regarding $\Diff(F\rel
F\cap T)$ as a subspace of $\Imb_f(F,\Sigma\rel T)$, we may assume
that the composition $\chi_3\chi_2$ is defined. On a sufficiently
small neighborhood of $i_P$ in $\Imb_t(P,\Sigma\rel T)$ define
$\chi(j)\in\Diff_f(\Sigma\rel T)$ by
$$\chi(j)=\chi_1(p(j))\,(\chi_3\chi_2)(\chi_1(p(j))^{-1}\circ j)\ .$$

\noindent Then for $x\in P$ we have
$$\eqalign{\chi(j)i_P(x)
&=\chi_1(p(j))\,(\chi_3\chi_2)(\chi_1(p(j))^{-1}\circ j)(x)\cr
&=\chi_1(p(j))\,\chi_1(p(j))^{-1} j (x)\cr
&=j(x)\cr}$$
\end{proof}

\noindent This yields immediately

\begin{corollary}{restrict imbeddings to S}{} Let $W$ be a vertical
suborbifold of $\Sigma$ containing $P$. Then $\Diff_f(\Sigma\rel
T)\to\Imb_t(P,\Sigma\rel T)$ and $\Imb_f(W,\Sigma\rel T)\to
\Imb_t(P,\Sigma\rel T)$ are locally trivial.
\end{corollary}

\section{The space of Seifert fiberings of a Haken 3-manifold}
\label{sfspace}
\marginwrite{sfspace}

Let $p\colon\Sigma\to {\cal O}$ be a Seifert fibering of a Haken
manifold $\Sigma$. As noted in section~\ref{sfiber}, $p$ is a singular
fibering.  Denote by $\diff_f(\Sigma)$ the connected component of the
identity in $\Diff_f(\Sigma)$, and similarly for other spaces of
diffeomorphisms and imbeddings. The main result of this section is the
following.

\begin{theorem}{space of fp homeos}{} Suppose that $\Sigma$ is a
Haken 3-manifold. Then the inclusion $\diff_f(\Sigma)\to
\diff(\Sigma)$ is a weak homotopy equivalence.
\marginwrite{space of fp homeos}
\end{theorem}

Before proving theorem~\ref{space of fp homeos}, we give an
application. Each element of $\Diff(\Sigma)$ carries the given
fibering to an isomorphic fibering, and $\Diff_f(\Sigma)$ is precisely
the stabilizer of the given fibering under this action. Therefore it
is reasonable to define the {\it space of Seifert fiberings}
isomorphic to the given fibering to be the space of cosets
$\Diff(\Sigma)/\Diff_f(\Sigma)$. Since $\Diff_f(\Sigma)$ is a closed
subgroup, the quotient $\Diff(\Sigma)\to
\Diff(\Sigma)/\Diff_f(\Sigma)$ is a principal fibering with fiber
$\Diff_f(\Sigma)$. As an immediate corollary to theorem~\ref{space of
fp homeos}, we will obtain:

\begin{theorem}{space of sf's}{} Suppose that $\Sigma$ is a
Haken 3-manifold. Then each path component of the space of Seifert
fiberings of $\Sigma$ is weakly contractible.
\marginwrite{space of sf's}
\end{theorem}

\begin{proof}{} As sketched on p.~85 of \cite{Waldhausen}, two
fiber-preserving diffeomorphisms of $\Sigma$ that are isotopic are
isotopic through fiber-preserving diffeomorphisms. This implies that
$\pi_0(\Diff_f(\Sigma))\to\pi_0(\Diff(\Sigma))$ is injective. By
theorem~\ref{space of fp homeos},
$\pi_q(\Diff_f(\Sigma))\to\pi_q(\Diff(\Sigma))$ is an isomorphism for
all $q\geq 1$. The theorem now follows from the homotopy exact
sequence for the fibration $\Diff(\Sigma)\to
\Diff(\Sigma)/\Diff_f(\Sigma)$.
\end{proof}

For compact Seifert fibered 3-manifolds, apart from a small list of
well-known exceptions, every diffeomorphism is isotopic to a
fiber-preserving diffeomorphism. So the following immediate corollary
applies to most cases.

\begin{corollary}{space of sf's coro}{} Suppose that $\Sigma$ is a
Haken 3-manifold such that every diffeomorphism is isotopic to
a fiber-preserving diffeomorphism. Then the space of Seifert fiberings
of $\Sigma$ is weakly contractible.
\marginwrite{space of sf's coro}
\end{corollary}

\noindent The proof of theorem~\ref{space of fp homeos} will use the
following lemma.

\begin{lemma}{rel fiber}{} Let $\Sigma$ be a  Haken
Seifert fibered 3-manifold, and let $C$ be a fiber of $\Sigma$. Then
each component of $\Diff_v(\Sigma\rel C)$ is contractible.
\marginwrite{rel fiber}
\end{lemma}

\begin{proof}{} Since $\Sigma$ is Haken, the base orbifold
of $\Sigma-C$ has nonpositive Euler characteristic and is not
closed. It follows (see~\cite{Scott}) that $\Sigma-C$ admits an
$\H^2\times\R$ geometry. Thus there is an action of $\pi_1(\Sigma-C)$
on $\H^2\times\R$ such that every element preserves the $\R$-fibers
and acts as an isometry on the $\H^2$ factor. Let $B$ be the orbit
space of~$\Sigma-C$.

It suffices to show that $\diff_v(\Sigma\rel C)$ is contractible. Let
$N$ be a fibered solid torus neighborhood of $C$ in $\Sigma$. It is
not difficult to see that $\diff_v(\Sigma\rel C)$ deformation retracts
to $\diff_v(\Sigma\rel N)$, which can be identified with
$\diff_v(\Sigma-C\rel N-C)$, so it suffices to show that the latter is
contractible. For $f\in\diff_v(\Sigma-C\rel N-C)$, let $F$ be a lift
of $f$ to $\H^2\times\R$ that has the form $F(x,s)\eq (x,s+F_2(x,s))$,
where $F_2(x,s)\in\R$. Since $f$ is vertically isotopic to the
identity relative to $N-C$, we may moreover choose $F$ so that
$F_2(x,s)\eq 0$ if $(x,s)$ projects to $N-C$. To see this, we choose
the lift $F$ to fix a point in the preimage $W$ of $N-C$. Since $f$ is
homotopic to the identity relative to $N-C$, $F$ is equivariantly
homotopic to a covering translation relative to $W$. That covering
translation fixes the point in $W$, and therefore must be the
identity. Thus $F$ fixes $W$ and commutes with every covering
translation.

Define $K_t$ by $K_t(x,s)\eq (x,s+(1-t)F_2(x,s))$. Since $K_0\eq F$
and $K_1$ is the identity, and each $K_t$ is the identity on the
preimage of $N-C$, this will define a contraction of
$\Diff_v(\Sigma-C\rel N-C)$ once we have shown that each $K_t$ is
equivariant. Let $\gamma\in\pi_1(\Sigma-C)$. From~\cite{Scott},
$\Isom(\H^2\times\R)=\Isom(\H^2)\times\Isom(\R)$, so we can write
$\gamma(x,s)=(\gamma_1(x),\epsilon_\gamma s+\gamma_2)$, where
$\epsilon_\gamma\eq \pm 1$ and $\gamma_2\in\R$. Since $F\gamma\eq
\gamma F$, a straightforward calculation shows that
$$F_2(\gamma_1(x),\epsilon_\gamma s+\gamma_2)\eq \epsilon_\gamma
F_2(x,s)\ .$$

\noindent Now we calculate
$$\eqalign{K_t\gamma(x,s)&= K_t(\gamma_1(x),\epsilon_\gamma
s+\gamma_2)\cr
&= (\gamma_1(x),\epsilon_\gamma s +\gamma_2+(1-t)%
F_2(\gamma_1(x),\epsilon_\gamma s+\gamma_2))\cr
&= (\gamma_1(x),\epsilon_\gamma s +\gamma_2+(1-t)%
\epsilon_\gamma F_2(x,s))\cr
&= (\gamma_1(x),\epsilon_\gamma (s+(1-t)F_2(x,s)) +\gamma_2)\cr
&= \gamma(x,s+(1-t)F_2(x,s))\cr
&= \gamma K_t(x,s)\cr}$$

\noindent showing that $K_t$ is equivariant.
\end{proof}

\begin{proof}{\ref{space of fp homeos}} We first examine
$\diff_v (\Sigma)$. Choose a regular fiber $C$ and consider the
restriction $\diff_v (\Sigma)\to \imb_v(C,\Sigma)\cong\diff(C)\cong
\diff(S^1)\simeq \SO(2)$. By corollary~\ref{sfcorollary2}(ii), this is
a fibration. By lemma~\ref{rel fiber}, each component of the fiber
$\Diff_v(\Sigma\rel C)\cap\diff_v(\Sigma)$ is contractible. It follows by the
exact sequence for this fibration that
$\pi_q(\diff_v(\Sigma))\cong\pi_q(\SO(2))\eq 0$ for $q\geq 2$, and for $q\eq
1$ we have an exact sequence $$0\longrightarrow
\pi_1(\diff_v(\Sigma))\longrightarrow
\pi_1(\diff(C))\longrightarrow
\pi_0(\Diff(\Sigma\rel C)\cap \diff_v(\Sigma))\longrightarrow
0\ .$$

We will first show that exactly one of the following holds.
\begin{enumerate}
\item[a)] $C$ is central and $\pi_1(\diff_v(\Sigma))\cong\Z$ generated by the
vertical $S^1$-action.
\item[b)] $C$ is not central and $\pi_1(\diff_v(\Sigma))$ is trivial.
\end{enumerate}

\noindent
Suppose first that the fiber $C$ is central in $\pi_1(\Sigma)$. Then
there is a vertical $S^1$-action on $\Sigma$ which moves the basepoint
(in $C$) once around $C$. This maps onto the generator of
$\pi_1(\diff(C))$, so $\pi_1(\diff_v(\Sigma))\to
\pi_1(\diff(C))$ is an isomorphism. Therefore
$\pi_1(\diff_v(\Sigma))$ is infinite cyclic, with generator
represented by the vertical $S^1$-action.

If the fiber is not central, then $\pi_1(\diff(C))\to
\pi_0(\Diff(\Sigma\rel C)\cap \diff_v(\Sigma))$ carries the
generator to a diffeomorphism of $\Sigma$ which induces an inner
automorphism of infinite order on $\pi_1(\Sigma,x_0)$, where $x_0$ is
a basepoint in $C$. Since elements of $\Diff(\Sigma\rel C)$ fix the
basepoint, this diffeomorphism (and its powers) are not in
$\diff(\Sigma\rel C)$. Therefore $\pi_1(\diff(C))\to
\pi_0(\Diff(\Sigma\rel C)\cap \diff_v(\Sigma))$ is injective, so
$\pi_1(\diff_v(\Sigma))$ is trivial.

Now consider the fibration of theorem~\ref{sfproject diffs}:
$$\Diff_v(\Sigma)\cap \diff_f(\Sigma)
\longrightarrow \diff_f(\Sigma)\longrightarrow
\diff({\cal O})\ .\leqno{(*)}$$

Observe that $\diff({\cal O})$ is homotopy equivalent to the identity
component of the space of diffeomorphisms of the 2-manifold ${\cal
O}-{\cal E}$, where ${\cal E}$ is the exceptional set. Since $\Sigma$
is Haken, this $2$-manifold is either a torus, annulus, disc with one
puncture, Mobius band, or Klein bottle, or a surface of negative Euler
characteristic. Therefore $\diff({\cal O})$ is contractible unless
$\chi({\cal O})=0$, in which case its higher homotopy groups are all
trivial, and its fundamental group is isomorphic to the center of
$\pi_1({\cal O})$. In the latter cases, the elements of $\pi_1({\cal
O})$ are classified by their traces at a basepoint of ${\cal O}-{\cal
E}$. From the exact sequence for the fibration $(*)$, it follows that
$\pi_q(\diff_f(\Sigma))\eq 0$ for~$q\geq 2$.

To complete the proof, we recall the result of Hatcher \cite{Hatcher}:
for $M$ Haken, $\pi_q(\diff(M))$ is $0$ for $q\geq 2$
and is isomorphic to the center of $\pi_1(M)$ for $q\eq 1$, and the
elements of $\pi_1(\diff(M))$ are classified by their traces at the
basepoint. We already have $\pi_q(\diff_f(\Sigma))\eq 0$ for
$q\geq 2$, so it remains to show that $\pi_1(\diff_f(\Sigma))\to
\pi_1(\diff(\Sigma))$ is an isomorphism.

\smallskip
\noindent {\sl Case I:} $\pi_1({\cal O})$ is centerless.
\smallskip

In this case $\diff({\cal O})$ is contractible, and either $C$
generates the center or $\pi_1(\Sigma)$ is centerless. The exact
sequence associated to the fibration $(*)$ shows that
$\pi_1(\diff_v(\Sigma))\to\pi_1(\diff_f(\Sigma))$ is an
isomorphism. Suppose $C$ generates the center. Since
$\pi_1(\diff_v(\Sigma))$ is infinite cyclic generated by the vertical
$S^1$-action, Hatcher's theorem shows that the composition
$$\pi_1(\diff_v(\Sigma)) \to
\pi_1(\diff_f(\Sigma))\to
\pi_1(\diff(\Sigma))$$

\noindent is an isomorphism. Therefore
$\pi_1(\diff_f(\Sigma))\to \pi_1(\diff(\Sigma))$ is an
isomorphism. If  $\pi_1(\Sigma)$ is
centerless, then $\pi_1(\diff(\Sigma))\eq 0$,
$\pi_1(\diff_f(\Sigma))\cong\pi_1(\diff_v(\Sigma))\eq 0$, and again
$\pi_1(\diff_f(\Sigma))\to \pi_1(\diff(\Sigma))$ is an
isomorphism.

\smallskip
\noindent {\sl Case II:} $\pi_1({\cal O})$ has center.

\smallskip
Assume first that ${\cal O}$ is a torus.
If $\Sigma$ is the $3$-torus, then by considering the exact sequence
for the fibration $(*)$, one can check directly that the homomorphism
$\partial\colon\pi_1(\diff({\cal
O}))\rightarrow\pi_0(\Diff_v(\Sigma)\cap\diff_f(\Sigma))$ is the zero
map. We obtain the exact sequence
$$0\longrightarrow \Z
\longrightarrow \pi_1(\diff_f(\Sigma))\longrightarrow
\Z\times\Z\longrightarrow0\ .$$
\noindent
Since $\diff_f(\Sigma)$ is a topological group,
$\pi_1(\diff_f(\Sigma))$ is abelian and hence isomorphic to
$\Z\times\Z\times\Z$. The traces of the generating
elements generate the center of $\pi_1(\Sigma)$, which shows that
$\pi_1(\diff_f(\Sigma))\to\pi_1(\diff(\Sigma))$ is an
isomorphism.

Suppose that $\Sigma$ is not a 3-torus. Then $\Sigma={\cal O}\times
I/(x,0)\simeq (\phi(x),1)$ for a homeomorphism $\phi\colon\cal O\rightarrow
\cal O$, $\pi_1(\Sigma)=\langle a,b,t \vbar tat^{-1}\eq a,
[a,b]\eq1,tbt^{-1}\eq a^nb\rangle$ for some integer $n$, and the fiber
$a$ generates the center of~$\pi_1(\Sigma)$.

Let $b_0$ and $t_0$ be the image of the generators of  $b$ and $t$
respectively in $\pi_1({\cal O})$. Now $\pi_1(\diff({\cal O}))\cong\Z\times\Z$
generated by elements whose traces represent the elements $b_0$ and $t_0$. By
lifting these isotopies we see that $\partial\colon\pi_1(\diff({\cal
O}))\rightarrow\pi_0(\diff_v(\Sigma))$ is injective. Therefore
$\pi_1(\diff_v(\Sigma))$ is isomorphic to  $\pi_1(\diff_f(\Sigma))$, and the
result follows as in case~I.

Assume now that ${\cal O}$ is a Klein bottle. As in the torus case
we may view $\Sigma={\cal O}\times I/(x,0)\simeq (\phi(x),1)$,
$\pi_1(\Sigma)=\langle a,b,t \vbar tat^{-1}\eq a^{-1},
[a,b]=1,tbt^{-1}=a^{-n}b^{-1} \rangle$ for some integer $n$, with
fiber $a$, and $\pi_1({\cal O})=\langle b_0,t_0 \vbar
t_0b_0t_0^{-1}\eq b_0^{-1}\rangle$. Now $\pi_1(\diff({\cal O}))$ is
generated by an isotopy whose trace represents the generator of the
center of $\pi_1(\diff({\cal O}))$, the element $t_0^2$. Observe that
$\pi_1(\Sigma)$ has center if and only if $n\eq 0$. If $n=0$, then it
follows that $\partial\colon\pi_1(\diff({\cal O}))\rightarrow
\pi_0(\Diff_v(\Sigma)\cap\diff_f(\Sigma))$ is the zero
map.  Hence $\pi_1(\diff_f(\Sigma))\rightarrow\pi_1(\diff(\cal O))$ is
an isomorphism and the generator of $\pi_1(\diff_f(\Sigma))$ is
represented by an isotopy whose trace represents the element $t^2$. By
Hatcher's result,
$\pi_1(\diff_f(\Sigma))\rightarrow\pi_1(\diff(\Sigma))$ is an
isomorphism. If $n\neq 0$, then $\partial\colon\pi_1(\diff({\cal
O}))\rightarrow\pi_0(\Diff_v(\Sigma)\cap\diff_f(\Sigma))$ is
injective. Since $\pi_1(\Sigma)$ is centerless,
$\pi_1(\Diff_v(\Sigma)\cap\diff_f(\Sigma))\eq0$. This implies that
$\pi_1(\diff_f(\Sigma))\eq0$, and again Hatcher's result applies.

The cases where $\cal O$ is an annulus, disc with one puncture, or a
Mobius band are similar to those of the torus and Klein bottle.
\end{proof}

\end{document}